\documentclass[reqno,12pt]{amsart}
\usepackage{amsfonts,a4}
\usepackage{amssymb ,bbold, bbm, marvosym}
\usepackage{latexsym}
\usepackage{amsmath}
\usepackage{amsthm}
\usepackage{color}
\usepackage{float}
\usepackage[mathscr]{eucal}
\usepackage{mathrsfs}
\usepackage[active]{srcltx}
\usepackage{tikz}
\usepackage{graphicx}
\font\sc=rsfs10 at 12pt

\numberwithin{equation}{section}

\newtheorem{theorem}{Theorem}[section]

\newtheorem{proposition}[theorem]{Proposition}
\newtheorem{corollary}[theorem]{Corollary}

\newtheorem{remark}[theorem]{Remark}

\newcommand {\B}    {\mathbb{B}}
\newcommand {\C}    {\mathbb{C}}
\newcommand {\Hb}   {\mathbb{H}}
\newcommand {\R}    {\mathbb{R}}

\newcommand {\Z}    {\mathbb{Z}}
\newcommand {\T}    {\mathbb{T}}
\newcommand{\Dn}     {\mathbb{D}}

\newcommand {\K}    {\mathbb{K}}

\newcommand{\rad}{\hbox{\footnotesize{\rm rad}\,}}

\renewcommand{\a}{\alpha}
\renewcommand{\b}{\beta}
\newcommand{\g}{\gamma}
\newcommand{\G}{\Gamma}

\newcommand{\D}{\Delta}
\newcommand{\e}{\epsilon}

\newcommand{\z}{\zeta}
\newcommand{\y}{\eta}

\renewcommand{\l}{\lambda}

\newcommand{\m}{\mu}

\renewcommand{\r}{\rho}

\renewcommand{\o}{\omega}

\renewcommand{\u}{\upsilon}
\newcommand{\by}{\pmb{\y}}
\newcommand{\kn}{\mathbb{k}}

\newcommand{\db}{{\mathbf d}}

\newcommand{\kb}{{\mathbf k}}

\newcommand{\nb}{{\mathbf n}}

\newcommand{\pb}{{\mathbf p}}

\newcommand{\sbb}{{\mathbf s}}

\newcommand{\ub}{{\mathbf u}}

\newcommand{\zb}{{\mathbf z}}

\newcommand{\Eb}{{\mathbf E}}

\renewcommand{\Hb}{{\mathbf H}}
\newcommand{\Ib}{{\mathbf I}}

\newcommand{\Pb}{{\mathbf P}}

\newcommand{\Tb}{{\mathbf T}}
\newcommand{\Ub}{{\mathbf U}}
\newcommand{\Vb}{{\mathbf V}}

\newcommand{\AF}{\mathfrak A}

\newcommand{\SF}{\mathfrak S}

\newcommand{\Ac}{{\mathcal A}}
\newcommand{\Bc}{{\mathcal B}}

\newcommand{\Dc}{{\mathcal D}}

\newcommand{\Hc}{{\mathcal H}}

\newcommand{\Rc}{{\mathcal R}}

\newcommand{\Tc}{{\mathcal T}}
\newcommand{\Vc}{{\mathcal V}} 

\font\sc=rsfs10 at 12pt

\newcommand{\Es}{\sc\mbox{E}\hspace{1.0pt}}

\newcommand{\Ls}{\sc\mbox{L}\hspace{1.0pt}}

\newcommand{\la}{\langle}
\newcommand{\ra}{\rangle}
\begin{document}

\date{\today}

\title[Commutative Toeplitz Algebras on the unit ball]{Spectral theorem approach to commutative $C^*$-algebras generated by Toeplitz operators on the unit ball: Quasi-elliptic related cases.}
\author{Grigori Rozenblum}
\address{Chalmers Univ. of Technol., Sweden; The Euler Intern. Math. Institute and St.Petersburg State Univ.; Mathematics Center
Sirius Univ. of Sci. and Technol.
Sochi Russia}
\email{grigori@chalmers.se}
\thanks{G.R. was supported  by the grant of the Russian Fund of Basic Research 20-01-00451.}
\author{Nikolai Vasilevski}
\address{Department of Mathematics, CINVESTAV, Mexico City, Mexico}
\email{nvasilev@math.cinvestav.mx}
\thanks{N.V. was partially supported by CONACYT grants 280732 and FORDECYT- \\
\indent PRONACES/61517/2020238630, Mexico.}

\begin{abstract}
We consider commutative $C^*$-algebras of Toeplitz operators in the weighted Bergman space on the unit ball in $\C^{\nb}.$ For the algebras of elliptic type we find a new representation, namely as the algebra of operators which are functions of  certain collections  of commuting unbounded self-adjoint operators in the Bergman space.

\vspace{1ex}
\noindent {\bf Keywords:} Toeplitz operators, Commutative algebras, Spectral representation
\vspace{1ex}\\
{\bf MSC (2020)}: Primary 30D60; Secondary 30G30; 30H20
\end{abstract}
\keywords Spectral functional calculus, Bergman spaces, Toeplitz operators
\maketitle

\section{Introduction}
An important topic in the analysis of Toeplitz operators in Bergman spaces is the description of commutative $C^*$- and Banach algebras generated by such operators. Such description was performed in the book \cite{Vas.book} and the series of papers, see e.g. \cite{GKV,GrudQuirogaVasil05,QuirogaVasil07a}, for weighted Bergman spaces on the unit disk and on the unit ball $\B^\nb$ in $\C^\nb$. Typical results in this direction relate commutative algebras with certain  groups of conformal transformations (biholomorphisms) of $\B^\nb$, $\nb \geq 1$, and the symbols of Toeplitz operators that generate commutative $C^*$-algebra, should be invariant with respect to these transformations.

Recently, in the paper \cite{RVSpectralDisk}, a different approach to the construction of commutative algebras was initiated. Namely, for a self-adjoint operator $V$ in the Bergman space  on the unit disk $\B^1$, the algebra of functions of this operator is commutative (in the sense corresponding to the size of the algebra of functions under consideration). If the operator $V$ commutes with a group of isometries of the Bergman space, generated by a one-parametric group of conformal mappings of the disk, the algebra of functions of this operator contains the $C^*$-algebra generated by Toeplitz operators, whose symbols are invariant under the above one-parametric group of conformal mappings. Such correspondence gives a better insight to the structure of  commutative algebras of operators,  extends the known classes of symbols generating bounded Toeplitz operators and, in particular, provides a clear explanation for the property of a Toeplitz operator to be compact.

In the present paper we extend the approach developed in \cite{RVSpectralDisk} to the study of commutative algebras of Toeplitz operators in the Bergman space on the ball in $\C^\nb.$ The passage to the multi-dimensional case requires some additional considerations. The most important is the one that the treatment of the classical commutative algebras of Toeplitz operators requires, generally,  not the usual spectral theory of self-adjoint operators but the multi-variable spectral theory for collections $\Vb=(V_1,\dots,V_m)$ of commuting self-adjoint operators.

Following the pattern of the one-dimensional case, it was shown in \cite{QuirogaVasil07a} that for each maximal Abelian subgroup (MASG) of the biholomorphisms of the unit ball $\mathbb{B}^{\nb}$, the $C^*$-algebra generated by Toeplitz operators whose symbols are invariant under the action of this subgroup, is commutative. Recall \cite{QuirogaVasil07a}, that there are five types of model MASGs on $\mathbb{B}^{\nb}$ (or its unbounded realization, the Siegel domain $\Dn_{\nb}$):

{\sf Quasi-elliptic} group of biholomorphisms of the unit ball $\mathbb{B}^{\nb}$, isomorphic to $\mathbb{T}^n$ with the  group action
\begin{equation*}
 t : \ z=(z_1,...,z_\nb) \in \mathbb{B}^\nb \longmapsto tz=(t_1z_1,...,t_\nb z_\nb) \in \mathbb{B}^\nb,
\end{equation*}
for each $t =(t_1,...,t_\nb )\in \mathbb{T}^\nb$.

{\sf Quasi-parabolic} group of biholomorphisms of the Siegel domain $\Dn_\nb$,  isomorphic to $\mathbb{T}^{\nb-1} \times \mathbb{R}$ with the  group action
\begin{equation*}
 (t,h): \ ( z',z_\nb) \in \Dn_\nb \longmapsto (tz',z_\nb+h) \in \Dn_\nb,
\end{equation*}
for each $(t,h) \in \mathbb{T}^{\nb-1} \times \mathbb{R}$.

{\sf Quasi-hyperbolic} group of biholomorphisms of the Siegel domain $\Dn_\nb$,  isomorphic to $\mathbb{T}^{\nb-1} \times \mathbb{R}_+$ with the  group action
\begin{equation*}
 (t,r) : \ (z',z_\nb) \in \Dn_\nb \longmapsto (r^{1/2}tz',rz_\nb) \in \Dn_\nb,
\end{equation*}
for each $(t,r) \in \mathbb{T}^{\nb-1} \times \mathbb{R}_+$,

{\sf Nilpotent} group of biholomorphisms of the Siegel domain $\Dn_\nb$,  isomorphic to $\mathbb{R}^{\nb-1} \times \mathbb{R}$
with the  group action
\begin{equation*}
 (b,h) : \ (z',z_\nb) \in \Dn_\nb \mapsto (z'+b,z_n+h+2iz'\cdot b +i|b|^2) \in \Dn_\nb,
\end{equation*}
for each $(b,h) \in \mathbb{R}^{\nb-1} \times \mathbb{R}$;

{\sf Quasi-nilpotent} group of biholomorphisms of the Siegel domain $\Dn_\nb$,  isomorphic to  $\mathbb{T}^{\kb} \times \mathbb{R}^{\nb-\kb-1} \times \mathbb{R}$, \ $0 < \kb < \nb-1$,
with the  group action
\begin{equation*}
 (t,b,h) : \ (z',z'',z_\nb) \in \Dn_\nb \longmapsto \\ (tz',z''+b,z_\nb+h+2iz'' \cdot b+i|b|^2)) \in \Dn_\nb,
\end{equation*}
for each $(t,b,h) \in \mathbb{T}^{\kb} \times \mathbb{R}^{\nb-\kb-1} \times \mathbb{R}$.

The latter, quasi-nilpotent, type of groups involves the integer parameter $\kb$, thus giving in total $\nb+2$ model groups for the $\nb$-dimensional ball $\mathbb{B}^{\nb}$.

In the present  paper we give a detailed characterization of commutative $C^*$-algebras generated by Toeplitz operators related to the \emph{quasi-elliptic group}  (and of  ambient commutative von Neumann algebras) as algebras of functions of the proper systems of commuting self-adjoint operators. The largest among the considered algebras is the algebra generated by Toeplitz operators with symbols being invariant under the smallest quasi-elliptic group $\mathbb{T}^{\nb}$, namely the diagonal subgroup of $\mathbf{U}(\nb)$ (the latter, as usual, denotes the group of all $\nb \times \nb$ unitary matrices). Here the collection of operators $\Vb$ consists of $\nb$ entries. Toeplitz operators with symbols invariant with respect to larger subgroups of $\mathbf{U}(\nb)$, the ones  consisting of block-diagonal matrices, generate commutative algebras coinciding with algebras of functions of smaller collections of commuting operators.  Such symbols are called \emph{quasi-radial}.

We consider also  Toeplitz operators with radial symbols, the ones invariant with the respect to the group $\mathbf{U}(\nb)$ itself. It turns out that this algebra is isomorphic to the algebra of functions of one single operator $\Vb=(V_1)$. Generally, the quantity of self-adjoint operators present in the functional calculus representation of the algebra generated by Toeplitz operators equals the number of diagonal  blocks in the subgroup of $\mathbf{U}(\nb).$

Finally, we consider from the spectral point of view certain commutative algebras of Toeplitz operators, which are proper subalgebras of the previous ones. Symbols of such Toeplitz operators  depend only on a part of variables, or possess some similar degeneracy. Such symbols and the corresponding Toeplitz operators were discussed, in particular, in \cite{BHaV,BV_2017,Quiroga_Sanchez_JFA}.  Contrary to the above, for each such algebra there is no subgroup of $\mathbf{U}(\nb)$ for which the symbols of the generating Toeplitz operators are exactly those which are invariant  under the action of this group. Nevertheless, a spectral representation of such algebras can be found.

\section{Preliminaries}
\subsection{Spectral Theorem}
For the  Readers' convenience  and in order to fix the notation, we recall some basic facts related to the Spectral Theorem for a finite collection of commuting self-adjoint operators. For proofs and more details see, e.g., \cite[Chapter 6, Section 5]{BirSol}, \cite[Chapter 1]{Samoilenko}, and \cite[Chapter 5]{Schmudgen}.

As  well known, for each, generally speaking unbounded, self-adjoint operator $V$ in a separable Hilbert space $\mathcal{H},$ there exists its spectral measure $E(\cdot)$, the projection-valued $\sigma$-additive function defined on  Borel sets in $\mathbb{R}$, orthogonal in the sense  $E(\Delta_1 \cap \Delta_2) = E(\Delta_1)E(\Delta_2)$,  and such that $E(\varnothing) = 0$, $E(\mathbb{R}) = \Ib$. Associated with $E(\Delta),$ is the resolution of identity, $\Eb(\eta)=E(-\infty,\eta].$ The operator $V$ admits the Stieltjes integral representation
\begin{equation*}
 V = \int_{\mathbb{R}} \eta d\Eb(\eta),
\end{equation*}
understood in the strong sense, with the  domain
\begin{equation*}
 \mathcal{D}(V) = \left\{f \in \mathcal{H}: \ \int_{\mathbb{R}} \eta^2 d\langle \Eb(\eta)f,f\rangle < \infty \right\}.
\end{equation*}
In the functional calculus associated with $V$, for a Borel complex-valued function $\psi(\eta)$ on $\R,$ the operator $\psi(V)$ is defined as
\begin{equation}\label{funct.1}
    \psi(V)=\int_{\R}\psi(\eta)d\Eb(\eta),
\end{equation}
again understood in the strong sense, with the domain
\begin{equation*}
    \Dc(\psi(V))=\left \{f \in \mathcal{H}: \ \int_{\mathbb{R}} |\psi(\eta)|^2 d\langle \Eb(\eta)f,f\rangle < \infty \right\}.
\end{equation*}

Such operators $\psi(V)$ are  normal, self-adjoint for a real-valued function $\psi$; they are bounded for $E$-essentially bounded functions $\psi$ and the mapping
\begin{equation}\label{Op.A}
    Op_{E}:\,\psi\mapsto \psi(V)
\end{equation}
is an isometric isomorphism of the commutative $C^*$ algebra of $E$-essentially  bounded Borel function to a commutative algebra of bounded operators.

For unbounded, $E$-almost everywhere finite functions $\psi,$ the normal operators $\psi(V)$ commute in the resolvent sense, as long as they possess at least one regular point.

Suppose now  that a finite collection of (possibly, unbounded) self-adjoint operators $\Vb=(V_1$, $V_2$, ..., $V_m)$ (an $m$-tuple) in the same separable Hilbert space $\mathcal{H}$ is given, and let $E_1(\cdot)$, $E_2(\cdot)$, ..., $E_m(\cdot)$ be their spectral measures.
We say that the operators $V_1$, $V_2$, ..., $V_m$ (strongly) commute if their spectral measures commute, i.e., for any Borel sets $\Delta_1$ and $\Delta_2$ in $\mathbb{R}$ and each $j'\neq j$, $E_j(\Delta_1) E_{j'}(\Delta_2)=E_{j'}(\Delta_2)E_j(\Delta_1)$.

\begin{theorem}[\sf Spectral Theorem]\label{Spectral Theorem}
For each $m$-tuple $\Vb = \{V_1,\ldots, V_m \}$ of pairwise stro\-ngly commuting self-adjoint operators in a separable Hilbert space  $\mathcal{H}$, there exists a unique \emph{joint spectral measure} $E(\cdot)$ on the $\sigma$-algebra of Borel sets in $\mathbb{R}^m$,  such that for the corresponding resolution of identity,
\begin{equation*}
 \Eb(\by)\equiv E\left(\prod_{j=1}^m (-\infty,\y_j]\right), \quad \by=(\y_1,\dots,\y_m),
\end{equation*}
the following Stieltjes representation holds:
\begin{equation*}
 V_j = \int_{\mathbb{R}^m} \eta_j d\Eb(\by), \quad j=1,\ldots,m,
\end{equation*}
with the integral understood in the strong sense.
\end{theorem}
The spectral measure $E(\cdot)$ can be constructed in the following way.
For a rectangle  $\Delta = [a_1,b_1] \times \ldots \times [a_m,b_m]$, we set $E(\Delta)= \prod_{j=1}^m E_j([a_j,b_j])$ and then extend it to all Borel sets on $\mathbb{R}^m$ in the standard way.

We recall also that the support of the spectral measure $E(\cdot)$ (the smallest closed set of full measure) is called the \emph{joint spectrum}
of the $m$-tuple $\Vb = \{V_1,\ldots, V_m \}$ and is denoted by
$\sigma(\Vb) \equiv \sigma(V_1,\ldots, V_m)\subset \R^m$; here
\begin{equation}\label{spectral product}
  \sigma(\Vb)=\prod_{j=1}^m\sigma(V_j).
\end{equation}

\begin{theorem}[\sf Functional Calculus]
Given an $E$-measurable function $\psi$ on $\mathbb{R}^m$, the operator
 \begin{eqnarray} \label{eq:calculus} \nonumber
 \psi(\Vb) &\equiv& \psi(V_1,\ldots, V_m) := \int_{\mathbb{R}^m} \psi(\eta_1,\ldots, \eta_m)\,d\Eb(\eta_1,\ldots, \eta_m) \\
 &=& \int_{\sigma(\Vb)} \psi(\eta_1,\ldots, \eta_m)\,d\Eb(\eta_1,\ldots, \eta_m)
 \end{eqnarray}
is well defined and \emph{normal} on its domain
\begin{equation*}
 \mathcal{D}_{\psi} = \left\{f \in \mathcal{H} : \int_{\mathbb{R}^m} |\psi(\eta_1,\ldots, \eta_m)|^2\, d\langle \Eb(\eta_1,\ldots, \eta_m)f,f\rangle < \infty \right\}.
\end{equation*}
The operator $\psi(V_1,\ldots, V_m)$ is bounded, and thus defined on the whole $\mathcal{H}$, if and only if the function   $\psi$ is $E$-essentially bounded.
The functional calculus mapping
\begin{equation*}
    Op_E: \psi\mapsto \psi(V_1,\ldots, V_m)\equiv\psi(\Vb)
\end{equation*}
is an isometric isomorphism of the $C^*$-algebra of $E$-essentially bounded complex functions onto a commutative $*$- algebra of of bounded normal operators in $\mathcal{H}.$
\end{theorem}
\begin{remark}\emph{In concrete situations, for a given tuple of operators $\Vb,$ the description of the range of $Op_E$ may be a serious analytical problem. This can be seen in the analysis in \cite{RVSpectralDisk} of the one-dimensional case}
\end{remark}

\subsection{The Bergman space on the ball} Our considerations deal with  objects in the complex space $\C^\nb$ for a certain $\nb.$ 
 Let $\mathbb{B}^\nb = \{ z=(z_1,...,z_\nb) \in \mathbb{C}^\nb : \ |z|^2 \equiv |z_1|^2 + ... + |z_\nb|^2 < 1\}$ be the unit ball in $\mathbb{C}^\nb$; for any $k\le \nb,$ for a ball $\mathbb{B}^k\subset \C^\nb,$ we denote by $\tau (\mathbb{B}^k)\subset \R^k\subset \C^\nb $ the base of the unit ball $\mathbb{B}^k$, considered as a Reinhard domain,  i.e.,
\begin{equation*}
    \tau(\mathbb{B}^k)= \{(r_1,..., r_k)=(|z_1|,...,|z_k|)\, : \
    r^2=r_1^2+...+r_k^2 \in [0,1) \}.
\end{equation*}
Given a multi-index $\alpha=(\alpha_1,\alpha_2,...,\alpha_\nb) \in \mathbb{Z}^\nb_+$  we will use the standard notation,
\begin{eqnarray*}
 |\alpha| &=& \alpha_1+\alpha_2+...+\alpha_\nb, \\
 \alpha! &=& \alpha_1!\,\alpha_2!\, \cdots\, \alpha_\nb!, \\
 z^{\alpha} &=& z_1^{\alpha_1}\,z_2^{\alpha_2}\,\cdots\,z_\nb^{\alpha_\nb},\\
 z_l=x_l+iy_l, \, l=1,2,...,\nb.
\end{eqnarray*}

Denote by $d\Vc=dx_1 dy_1 ...dx_\nb dy_\nb$ the standard Lebesgue measure in $\mathbb{C}^\nb$; for $\lambda > -1$, we introduce the probability weighted measure on $\B^{\nb},$
\begin{equation*}
dv_\lambda(z) = c_{\lambda}\, (1-|z|^2)^{\lambda}\,d\Vc(z),
 \quad \mathrm{where} \quad c_{\lambda} = \frac{\Gamma(\nb+\lambda+1)}{\pi^\nb\, \Gamma(\lambda+1)}, \ \ \ \lambda > -1.
\end{equation*}
The basic object of our considerations is  the weighted Hilbert space $\Hc^\nb_\l\equiv L_2(\mathbb{B}^\nb, dv_{\lambda})$ and its subspace, the weighted Bergman space $\mathcal{A}^2_{\lambda}=\mathcal{A}^2_{\lambda}(\mathbb{B}^\nb)$ which consists of all analytic functions  in~$\Hc^\nb_\l$. By $\mathbb{P}$ we denote the orthogonal projection onto $\mathcal{A}^2_{\lambda}$.

Further on, in notations,  the superscript $2$ will be skipped since our considerations concern  only the Hilbert space theory.

Recall that the standard orthonormal basis of the Bergman space $\Ac_\l(\B^{\nb})$ consists of  the monomials
\begin{equation}\label{basis}
e_{\alpha}(z):= e_{\alpha,(\nb)}(z) \equiv \o_{\l,\a,(\nb)}z^{\alpha}, \quad \mathrm{with} \quad \alpha \in \mathbb{Z}_+^\nb.
\end{equation}
where
$\o_{\l,\a,(\nb)}= \sqrt{\frac{\Gamma(\nb+|\alpha|+\lambda+1)}{\alpha! \Gamma(\nb+\lambda+1)}}.$
 Note that the normalization factor $\o_{\l,\a,(\nb)}$ in \eqref{basis} depends on the dimension $\nb$, so the function $z^\a$ having $\a_l=0$ for $l>k$ with some $k$, $1\le k<\nb$, produces an element $e_{\alpha,(k)}(z)$ in the basis of the space $\Ac_\l(\mathbb{B}^k)$, but with the normalization factor $\o_{\l,\a,(k)},$ different from the one in \eqref{basis}. If some value $k$ different from $\nb$ comes into consideration, the subscript $(k)$ will be present in the notation $e_{\a,(k)}$, otherwise, $k=\nb$ is assumed. In such sections where the weight exponent $\l$ is fixed, we may  omit it in the notation as well.


Thus, the space $\Ac_\l(\B^\nb)=\Ac(\B^\nb)$ consists of analytic functions $f(z)=f(z_1,\dots,z_\nb)$ for which the Fourier coefficients
\begin{equation}\label{Fourier}
    f_\a=\langle f,e_{\alpha}\rangle_{\Hc_{\l}(\B^\nb)}
\end{equation}
satisfy $\sum_{\a\in\Z^\nb_+}|f_\a|^2<\infty.$

\subsection{Partitions and Isometry Groups}\label{isometry} 


Let $\kn=(k_1,...,k_m)$ be a \emph{tuple} of positive integers: $k_1+...+k_m =\nb$.
Given the partition $\kn=(k_1,...,k_m)$, we arrange the coordinates  of $z \in \mathbb{B}^\nb$ in $m$ groups, each one of which has $k_j$, $j=1,...,m$, entries and introduce the notation
\begin{equation*}
 z_{(1)}= (z_{1,1},...,z_{1,k_1})\in\C^{k_1}, \ z_{(2)}= (z_{2,1},...,z_{2,k_2})\in \C^{k_2}, \
..., \ z_{(m)}= (z_{m,1},...,z_{m,k_m})\in\C^{k_m},
\end{equation*}
with a proper re-numbering of co-ordinates, so that
\begin{equation*}
 \mathbb{B}^\nb \ni z = (z_1,z_2,\ldots,z_\nb) = (z_{(1)}, z_{(2)}, \ldots, z_{(m)}).
\end{equation*}
With a group $z_{(j)}$ fixed, the \emph{complementing} group of variables will be denoted by
\begin{equation}\label{z-complement}
    \widetilde{z_{(j)}}=(z_{(1)},\ldots,z_{(j-1)},z_{(j-1)},\ldots,z_{(m)}),
\end{equation}
so, with an obvious permutations of variables (which is always assumed performed), $z=(z_{(j)}, \widetilde{z_{(j)}}).$ In the extreme case $\kn=(1,1,\dots,1),$ $\widetilde{z_j}=\widetilde{z_{(j)}}.$
(Note the subtlety in notation: the subscript $j$ (as in $z_j$, $\widetilde{z_j}$, $V_j$, etc.) is used for denoting a single complex variable or some related object, while the subscript $(j)$ is used to denote a tuple of variables and related objects, as in $z_{(j)},\widetilde{z_{(j)}} $ etc.)

 With a partition  $\kn=(k_1,...,k_m)$ we associate the subgroup $\Hb=\Hb_{\kn}\subset \Ub(\nb)$
\begin{equation*}
 \Hb_{\kn} = \mathbf{U}(\kn) := \mathbf{U}(k_1) \times \ldots
 \times \mathbf{U}(k_m),
\end{equation*}
realized as a block diagonal subgroup of $\mathbf{U}(\nb)$.

Among such subgroups, the smallest one and the only commutative one
is the group $\Hb_{\min}$ corresponding to the partition $\kn_{\min}=(1,1,\dots,1)$; this group is the $\nb$-dimensional torus $
\Hb_{\min}=\T^{\nb}$ realised as the group of diagonal unitary matrices. The largest of such subgroups is $\Ub(\nb)$ itself. It corresponds to the trivial partition $\kn=(\nb).$ All other partitions generate groups partially ordered by inclusion.

With each partition $\kn$ we associate the class of (first, bounded) symbols $\SF(\kn)$ invariant under the action of the group $\Hb_{\kn}.$  It stands to reason that the classes  $\SF(\kn)$ are ordered in the inverse way to the groups $\Hb_{\kn}$. As we will see, Toeplitz operators with symbols  in $\SF(\kn)$ generate a commutative $C^*$-algebra. In the sections to follow  we associate with each subgroup $\Hb_{\kn}$ a collection of commuting unbounded self-adjoint operators $\Vb=\Vb(\kn)$ so that the (von Neumann) algebra of all bounded measurable functions of the operators $\Vb$ is commutative and contains the $C^*$-algebra generated by Toeplitz operators with symbols in $\SF(\kn)$.
It turns out that the smaller is the group $\Hb_{\kn},$ the more operators in the system $\Vb$ are needed. In the extreme cases $\Hb=\Ub(\nb),$ respectively, $\Hb=\T^{\nb}$, the system $\Vb$ consists of one, respectively, $\nb$ operators.

\section{Toeplitz operators with separately radial symbols} We consider first the subgroup $\Hb_{\min}=\T^{\nb}$ represented as the group of diagonal unitary  matrices in a certain fixed co-ordinate system $\zb=(z_1,\dots,z_{\nb})$ in $\B^{\nb}.$ This group, the smallest one, corresponds to the maximal partition, $\kn_{\max}=(1,\dots,1).$ Symbols in the class  $\SF(\kn_{\max})$ are invariant with respect to the action of $\T^{\nb};$ these are functions $a(r_1,\dots,r_\nb),$ where $r_j=|z_j|,$ $j=1,\dots, \nb.$ Such symbols are called \emph{separately radial}.

Following our convention, for $z=(z_1,\dots,z_{\nb})\in \B^{\nb}$ and a fixed $j\in[1,\nb]$ we denote by $\widetilde{z_{j}}$ the complementing collection of variables,
\begin{equation*}
   \widetilde{z_{j}}=(z_1,\dots, z_{j-1},z_{j+1},\dots z_{\nb}),
\end{equation*}
thus, with an obvious permutation, $z=(z_j, \widetilde{z_j}).$

 A basic  function $e_\a=e_{\a,(\nb)}$ is not separately radial, but the one-dimensional subspace spanned by $e_\a$ \emph{is} invariant with respect to the action of $\T^{\nb}$, and this property pertains after the multiplication by a separately radial function $a(z)\in\SF(\kn_{\max})$, as well as the orthogonality relations,
\begin{equation*}
     \la e_\a,a(z)e_{\a'}\ra=0, \, \a\ne\a', a\in\SF(\kn_{max}).
\end{equation*}
Therefore, for a separately radial symbol $a$, the Toeplitz operator $\Tb_a: f\mapsto \mathbb{P}af$ in $\Ac_\l(\B^\nb)$ is diagonal in the basis $e_\a,$
 \begin{equation}\label{Separ.basis}
    \Tb_a e_\a =\g_{a}e_\a=\g_{a,\mathrm{sep}}e_\a,
 \end{equation}
 where
 \begin{gather}\label{coeff.separately}
 \g_{a}(\a)\equiv\g_{a,\mathrm{sep}}(\a)=\\\nonumber
 \ \frac{2^\nb\, \Gamma(\nb+|\alpha| + \lambda + 1)}{\Gamma(\lambda+1)  \a!}
 \int_{\tau(\mathbb{B}^\nb)} a(r_1,...,r_\nb) (1-|r|^2)^{\lambda}  \prod^\nb_{j=1} r_j^{2|\alpha_{j}|+1} dr_j \\ \nonumber
 = \ \frac{ \Gamma(\nb+|\alpha| + \lambda + 1)}{\Gamma(\lambda+1) \a!}
 \int_{\pmb{\Delta}_{\nb}} a(\sqrt{\r_1},...,\sqrt{\r_\nb}) (1-(\r_1+\ldots + \r_\nb))^{\lambda}  \\ \nonumber
 \times \ \prod^\nb_{j=1} \r_j^{\alpha_{j}} d\r_j;
 \end{gather}
 here $\pmb{\Delta}_{\nb}=\{\r\in\R^{\nb}_+:\r_1+\dots+\r_{\nb}<1\}$ is the standard $\nb$-dimensional simplex, $\r_j=r^2_j$.

 For a fixed multi-index $\sbb=(s_1,\dots,s_{\nb})\in\Z^\nb_+$, we define, following ( \cite[Formula (3.10)]{BV_2013}), the one-dimensional subspace $H_\sbb$ spanned by the basic function $e_{\sbb}.$ We also define, for a given $\ell\in \Z_+$ and a fixed $j,$ the \emph{infinite-dimensional} subspace $H_\ell^{(j)}\subset \Ac_\l(\B^\nb)$ as
 \begin{equation}\label{Compl H ell}
    H_\ell^{(j)}\equiv\overline{\check{H}_\ell^{(j)}} \equiv\overline{ \mathrm{span}}\,\{e_{\alpha}: \ \alpha_j =\ell \}.
 \end{equation}
 It is important to exlpain here that the non-closed space ${\check{H}_\ell}^{(j)}$ in \eqref{Compl H ell} consists of all polynomials in $z$ containing the variable $z_j$ taken exactly   to the power $\ell,$ $f(z)=z_j^{\ell}h(\widetilde{z_{j}}).$  After being  closed in the norm of $\Ac_\l(\B^{\nb})$, this structure is modified in the following way. We calculate the norm in  $\Ac_\l(\B^{\nb})$ of a function  $f(z)=z_j^{\ell}h(\widetilde{z_{j}})$.
 We have
 \begin{gather}\label{calculation1}
    \|f\|^2_{\Ac_\l(\B^{\nb})}=c_\l\int_{\B^\nb}|z_j^{\ell}|^2|h(\widetilde{z_{j}})|^2 (1-|z_j|^2-|\widetilde{z_{j}}|^2)^{\l}d\Vc(z)\\
    \nonumber
   = c_\l \int_{\B^{\nb-1}}|h(\widetilde{z_{j}})|^2 d\Vc(\widetilde{z_{j}}) \int_{|z_j|<(1-|\widetilde{z_{j}}|^2)^{\frac12}}|z_j|^{2\ell}(1-|z_j|^2-|\widetilde{z_{j}}|^2)^{\l} d \Vc(z_j).
 \end{gather}
We set $|\widetilde{z_{j}}|=\widetilde{r_{j}}$ and calculate the $z_j$-integral in \eqref{calculation1}, making use of \cite[Formula 3.191.1]{GR}
\begin{gather}\label{calculation2}
    \int_{|z_j|<(1-\widetilde{r_{j}}^2|)^{\frac12}}|z_j|^{2\ell}(1-|z_j|^2-|\widetilde{z_{j}}|^2)^{\l}d\Vc(z_j)\\ \nonumber
  =  2\pi \int_{0}^{(1-\widetilde{r_{j}}^2)^{\frac12}}r_j^{2\ell+1}(1-r_j^2-\widetilde{r_{j}}^2)^{\l}dr_j\\ \nonumber
   = \pi\int_{0}^{1-\widetilde{r_{j}}^2}\r_j^{\ell}(1-\r_j-\widetilde{r_{j}}^2)^\l d\r_j \\ \nonumber
  =  \pi (1-|\widetilde{r_{j}}^2|)^{\ell+\l+1} \frac{\G(\ell+1)\G(\l+1)}{\G(\ell+\l+2)}.
\end{gather}
Thus, the square of the norm of $f=z_j^{\ell}h(\widetilde{z_{j}})$ in $\Ac_\l(\B^{\nb})$ equals
\begin{equation}\label{normchange}
    \|f\|^2_{\Ac_\l(\B^{\nb})}=c_{\l,\ell}\|h\|_{{\Ac_{\l+\ell+1}(\B^{\nb-1})}}^2,
\end{equation}
with
\begin{equation}\label{coeff.change}
    c_{\l,\ell} = (\nb+\lambda)\frac{\G(\ell+1)\G(\l+1)}{\G(\ell+\l+2)}
\end{equation}

As a result of these calculations we can see that the space  $H_\ell^{(j)}$ defined in \eqref{Compl H ell} consists of functions of the form $f(z)=z_j^\ell h(\widetilde{z_{j}}),$ with $h$ belonging to the Bergman space $\Ac_{\l+\ell+1}(\B^{\nb-1})$. Moreover, the norms of $f$ in $\Ac_\l(\B^{\nb})$ and of $h$ in $\Ac_{\l+\ell+1}(\B^{\nb-1})$ are equivalent.

 It follows from the definition that $H_\sbb=\bigcap_{j=1}^{\nb}H_{s_j}^{(j)}$ for $\sbb\in\Z^\nb_+.$
We introduce the orthogonal projections, $\Pb_{\sbb}:\Hc\longrightarrow H_{\sbb}$ (a rank one projection) and $\Pb_\ell^{(j)}: \Hc\longrightarrow H_{\ell}^{(j)},$  so that $\Pb_\sbb=\prod_{j=1}^{\nb} \Pb_{s_j}^{(j)}.$
For each $j,\ell,$ the monomials $e_{\a}$ with $\a_j=\ell$ form an orthonormal basis in  $H_\ell^{(j)}$. Therefore the projection $\Pb_\ell^{(j)}$ can be written as
 \begin{equation}\label{P_l^j}
    (\Pb_\ell^{(j)}f)(z_j,\widetilde{z_{j}})=\sum_{\a_j=\ell}\la f,e_\a\ra e_{\a}(z).
\end{equation}
 With this notation, the Bergman space $\Ac_\l(\B^\nb)$ splits into the  orthogonal sum in two ways:
 \begin{equation}\label{OrtSumSepar}
    \Ac_\l(\B^{\nb})=\bigoplus_{\sbb\in\Z^{\nb}_+} H_{\sbb} \qquad  \mbox{and} \qquad \mbox{for each}\, j=1,\dots, \nb,\qquad \Ac_\l(\B^{\nb})=\bigoplus_{\ell\in\Z_{+}} H_{\ell}^{(j)}.
 \end{equation}
 We can see that in the first decomposition, we split the Bergman space into the direct sum of one-dimensional subspaces,  in the second one each subspace in the splitting consists of functions having a fixed homogeneity order  in the variable $z_j$ (and analytic in all remaining variables.)

For $f$ being a basic  monomial, $f(z)=e_\a(z_j,\widetilde{z_{j}})$, the projection $\Pb_\ell^{(j)}$ acts, according to \eqref{P_l^j},
as the identity operator in $\widetilde{z_{j}}$ variables and as an integral operator in $z_j$ variable. By continuity, this description extends to the whole $H^{(j)}_{\ell}.$ Therefore, for $j'\ne j$ and any $\ell',\ell\in \Z_+$ the projections $\Pb_\ell^{(j)}$ and $\Pb_{\ell'}^{(j')}$ commute; taken in any order, their composition is the projection onto the space of functions of the form $z_j^{\ell}z_{j'}^{\ell'}h(\widetilde{z_{j,j'}}),$ where $\widetilde{z_{j,j'}}$ is the collection of variables complementary to $z_j,z_{j'}$ and $h$ is a function in $\Ac_{\l+\ell+\ell'+2}(\B^{\nb-2}).$

 Now we introduce the \emph{rotation operator} with respect to the variable $z_j,$
 \begin{equation}\label{Rotation.Separate}
    V_j=z_j\frac{\partial}{\partial z_j}=\frac{1}{\imath}\frac{\partial}{\partial \theta_j},
 \end{equation}
 in the trigonometrical representation $z_j=|z_j|e^{\imath \theta_j}.$
 This operator is self-adjoint being considered on the domain
 \begin{equation*}
    \Dc(V_j)=\left\{f=\sum_{\a\in\Z^{\nb}_+}f_{\a}e_{\alpha} \in \Ac_\l(\B^\nb): \sum_{\a\in\Z^{\nb}_+}|\a_j|^2|f_\a|^2<\infty\right\}.
 \end{equation*}
For a fixed $j$, each element in the space $H_\ell^{(j)}$ is an eigenfunction of $V_j$ with eigenvalue $\ell$. The second orthogonal sum decomposition in \eqref{OrtSumSepar} implies that  finite linear combinations of elements in all $H_\ell^{(j)}$ are dense in $\Ac_\l(\B^\nb)$.

 For any fixed $j,$ each function of the form $z_j^\ell h(\widetilde{z_{(j)}}),$ $h\in\Ac_{\l+\ell+1}(\B^{\nb-1}),$ is an eigenfunction of $V_j$ with eigenvalue $\ell;$ linear combinations of all such functions are dense in $\Ac_\l(\B^\nb)$. This means that the spectrum of $V_j$ in $\Ac_\l(\B^\nb)$ consists of all nonnegative integer points $\ell\in \Z_+$ with the corresponding spectral subspace $H_\ell^{(j)};$ each integer  point is an eigenvalue of infinite multiplicity.

It follows that the spectral measures of the self-adjoint  operators $V_j,$ $j\le \nb$ commute, with $\Z_+^{\nb}$ being the joint spectrum. Therefore, the joint  spectral measure $E$ is supported on the integer lattice $\Z^\nb_+$, while the set $\R^\nb \setminus \Z^\nb_+$ has  zero $E$-measure. The corresponding partition of unity is
\begin{equation*}
  \Eb(\y_1,\dots,\y_\nb)=\prod_{j=1}^\nb\left(\sum_{s_j\le \y_j}\Pb_{s_j}^{(j)}\right).
\end{equation*}

 The general spectral theory for systems of commuting operators, presented in Sect.2, applies to $\Vb=(V_1,\dots,V_{\nb})$. 

Therefore, the Functional Calculus produces
\begin{proposition}
 Given an $E$-measurable function $\psi(\by)$ on $\mathbb{R}^\nb$, the operator
 \begin{equation*}
  \psi(\Vb) \equiv\psi(V_1,\ldots, V_\nb) := \int_{\mathbb{R}^\nb} \psi(\eta_1,\ldots, \eta_\nb)\,d\Eb(\eta_1,\ldots, \eta_\nb)
 \end{equation*}
is well defined and normal on its domain
\begin{equation*}
 \mathcal{D}_{\psi} = \left\{f \in \mathcal{\Ac_{\l}(\B^{\nb})} : \int_{\mathbb{R}^\nb} |\psi(\eta_1,\ldots, \eta_\nb)|^2\, d\langle \Eb(\eta_1,\ldots, \eta_\nb)f,f\rangle < \infty \right\}.
\end{equation*}
Moreover the operator $\psi(V_1,\ldots, V_\nb)$ is bounded, and thus defined on the whole $\mathcal{H}$, if and only if the function   $\psi$ is $E$-essentially bounded.
\end{proposition}

Representation \eqref{eq:calculus} implies that
\begin{equation*}
 \psi(\mathbf{V}) \equiv \psi(V_1,\ldots, V_\nb) = \int_{\sigma(\mathbf{V})} \psi(\eta_1,\ldots, \eta_n)\,d\Eb(\eta_1,\ldots, \eta_\nb) = \sum_{\sbb \in \Z^\nb_+} \psi(\sbb)\Pb_\sbb.
\end{equation*}

\begin{corollary} \label{co:R}
 The set of all operators $\{\psi(\mathbf{V})\}$, defined by (the classes of equivalence) of  $E$-measurable essentially bounded functions  $\psi$ with $\boldsymbol\psi= \{\psi(\sbb)\}_ {\sbb \in \mathbb{Z}^\nb_+} \in \ell_{\infty}(\mathbb{Z}^\nb_+)$ constitutes an algebra $\mathcal{R}$ of bounded mutually commuting operators in $\mathcal{A}_{\lambda}(\mathbb{B}^{\nb})$, and the mapping
 \begin{equation*}
  Op_{E}:\,\boldsymbol\psi \  \longmapsto \psi(\mathbf{V}) = \sum_{\sbb \in \Z^\nb_+} \psi(\sbb)P_\sbb.
 \end{equation*}
defines the isomorphism
\begin{equation}\label{OP,E}
 Op_{E}:\, \ell_{\infty}(\mathbb{Z}^\nb_+) \ \longrightarrow \ \mathcal{R}
\end{equation}
of the $C^*$-algebras.
\end{corollary}

Note here that although the operator $V_j$ acts upon the basis elements in the same way  in weighted spaces  $\Ac_\l(\B^{\nb})$ for all values of $\l>-1,$ the domain of $V_j$ depends on $\l,$ therefore, the self-adjoint operators $V_j$ as well as their spectral projections depend on $\l,$ although this dependence might be not reflected in the notation.

Returning now to Toeplitz operators with separately radial symbols, we see that by \eqref{Separ.basis}, the Toeplitz operator $\Tb_a$ with such $a= a(r_1,r_2,\ldots,r_\nb)$ admits the representation
\begin{equation*}
 T_a = \sum_{\sbb \in \mathbb{Z}_+^\nb} \gamma_a(\sbb) P_\sbb,
\end{equation*}
where the scalars $\gamma_a(\sbb)$ are given by \eqref{coeff.separately}. Thus, by Corollary \ref{co:R},
\begin{equation*}
 \Tb_a = \psi_a(V_1,V_2,\ldots,V_\nb),
\end{equation*}
where $\psi_a$ belongs to the equivalence class of essentially bounded $E$-measurable functions, defined by the sequence $\boldsymbol\psi_a= \{\psi_a(\sbb)\}_ {\sbb \in \mathbb{Z}^\nb_+}$, with $\psi_a(\sbb) = \gamma_a(\sbb)$.
In other words
\begin{equation}\label{Toeplitz.Decomp.Separ}
    \Tb_a=\sum_{\sbb\in\Z^\nb_+}\psi_a(\sbb)\Pb_{\sbb} = \sum_{\sbb\in\Z^\nb_+}\gamma_a(\sbb)\Pb_{\sbb}.
 \end{equation}

\begin{theorem}\label{sep.th} Denote by $\Tc(\kn_{\max})$ the $C^*$-algebra generated by  Toeplitz operators with bounded separately radial symbols. Then $\Tc(\kn_{\max})$ is contained in the image $\Rc$ in the mapping \eqref{OP,E}. In other words, every operator in $\Tc(\kn_{\max})$ is a certain function of the system of commuting operators $V_j.$
\end{theorem}

Similar to the one-dimensional case, $\Tc(\kn_{\max})$ does not exhaust $\Rc,$ this means that there exist  functions $\psi$ such
 that the operator $\psi(\Vb)$ does not belong to the $C^*$ algebra $\Tc(\kn_{\max})$, although $\psi$ is bounded. On the other hand, a Toeplitz operator with an unbounded separately radial symbol $a$ can belong to $\Rc.$ Corresponding examples can be constructed  on the base of reasoning in \cite{RVSpectralDisk}. Note also that Theorem \ref{sep.th} implies a compactness condition for Toeplitz operators of the type under consideration. Namely, sice all points of spectrum of $\Vb$ are eigenvalues of finite multiplicity placed at the nods of a lattice, an operator $\psi(\Vb)$ is compact iff $\psi(\sbb)\to 0$ as $|\sbb|\to\infty.$ For a Toeplitz operator, such behavior of $\sbb$ can be derived from the convergence to zero in a proper sense of the symbol, as $|z|\to 1,$ similarly to \cite{RVSpectralDisk}.  In Sect.6 we encounter the situation when this kind of results on the compactness break down.


\section{Toeplitz operators with radial symbols}\label{Sect.radial}
 The radial symbols correspond to the smallest partition $\kn_{min}=(\nb)$ consisting of only one element. The corresponding isometry group is $\Ub(\nb)$. Although this group is not commutative,  Toeplitz operators $\Tb_{a}$ with symbols in $\SF(\kn_{\min})$, commute, see \cite{Vas.book}. Such symbols, invariant with respect to $\Ub(\nb)$, depend only on the radius $r=(\sum r_j^2)^{\frac12}$ and therefore they are called 'radial'. We denote by $\AF_{\rad}$ the $C^*$ algebra generated  by Toeplitz operators with bounded symbols in $\SF(\kn_{\min})$. An extensive description of properties of operators with radial symbols can be found in \cite{BHV}. We use these facts to perform the spectral analysis of such operators, following the lines of our paper \cite{RVSpectralDisk}, where the case of the Bergman space on the disk $(\nb=1)$ was considered. We will see that, similarly to the case of the disk, the spectral representation of the algebra $\AF$ uses only one operator $V$.

 We define the differential operator (the rotation operator)
\begin{equation}\label{Operator V}
    V=z_1\frac{\partial}{\partial z_1}+\dots +z_{\nb}\frac{\partial}{\partial z_\nb}=\frac{1}{\imath}\frac{\partial}{\partial \theta_1}+\dots+\frac{1}{\imath}\frac{\partial}{\partial \theta_\nb},
\end{equation}
in the trigonometric representation $z_j=|z_j|e^{\imath\theta_j}$.
This operator is self-adjoint in $\Ac_{\l}(\mathbb{B}^\nb),$ being defined on
\begin{equation}\label{D(V)}
    \Dc(V)=\left\{f=\sum_{\a\in\Z_+^{\nb}}f_\a e_{\a}\in\Ac_{\l}(\B^{\nb}):\sum_{\a\in\Z_+^{\nb}}|\a|^2|f_\a|^2<\infty\right\}.
\end{equation}

The operator $V$ is diagonal in the basis $\{e_{\a}\}$ and it acts on an element $e_\a$ in the basis as $Ve_{\a}=|\a|e_{\a}.$ Therefore, each nonnegative integer $\ell\in \Z_+$ is an eigenvalue of $V,$ and the multiplicity of $\ell$ equals the number of multiindices $\a\in\Z_+^\nb$ satisfying $|\a|=\ell$. As well known, this number equals $\db_\ell=\binom{\ell+\nb-1}{\nb-1}.$
We denote by $\Ls_\ell$ the eigenspace of $V$ corresponding to the eigenvalue $\ell$, so $\dim\Ls_\ell=\db_\ell.$

Let $\Pb_\ell$ be the orthogonal (in $\Ac_{\l}(\B^\nb)$) projection onto $\Ls_\ell.$ It can be expressed as
\begin{equation}\label{P_l}
   \Pb_\ell f=\sum_{|\a|=\ell} \la f,e_{\a}\ra e_{\a},
\end{equation}
where $\la .,.\ra$ denotes the scalar product in the Hilbert space $\Hc_\l(\B^{\nb}).$ By \eqref{P_l}, $\Pb_\ell$ is an integral operator
\begin{equation*}
    (\Pb_\ell f)(z)=\int_{\B^\nb}f(\z)\Pb_\ell(z,\z)dv_\l(\z)
\end{equation*}
with the (degenerate) integral kernel
\begin{equation}\label{P_l kernel}
    \Pb_\ell(z,\z)=\sum_{|\a|=\ell}e_{\a}(z)\overline{e_{\a}(\z)}=\sum_{|\a|=\ell}\pb_\a z^\a\bar{\z}^\a
\end{equation}
where $\pb_\a\equiv\pb_{\a,(\nb)}=\frac{1}{\a!}\frac{\G(\nb+\ell+\l+1)}{\G(\nb+\l+1)},$
in accordance with \eqref{basis}.
Since
\begin{gather*}
    \sum_{|\a|=\ell}\frac{1}{\a!}z^\a\bar{\z}^\a=\frac{1}{\ell!}\sum_{|\a|=\ell}\frac{\ell!}{\a!}z^\a\bar{\z}^\a=\\\nonumber
  \frac{1}{\ell!}\sum_{|\a|=\ell} \frac{\ell!}{\a!} \prod_{1\le j \le \nb}(z_j\bar{\z_j })^{\a_j}=\frac{1}{\ell!}\la z,\bar{\z}\ra^\ell,
\end{gather*}
\eqref{P_l kernel} gives
\begin{equation}\label{P_l kernl fin}
    \Pb_\ell(z,\z)=\frac{1}{\ell!}\frac{\G(\nb+\ell+\l+1)}{\G(\nb+\l+1)}\la z,\bar{\z}\ra^\ell.
\end{equation}

The spectral measure for $V$ is given by
\begin{equation}\label{Sp.Meas.V}
    E_V(\D)=\sum_{\ell\in\D}\Pb_\ell,
\end{equation}
with the corresponding resolution of identity
\begin{equation}\label{Res.Id.V}
    \Eb_V(\eta) = \sum_{\ell\leq\eta}\Pb_\ell.
\end{equation}
Thus, for any $f\in\Ac_\l(\B^\nb)$,
the measure $\m_f(\D)=\la E_V(\D) f, f\ra$ is supported at integer points.

Further on,  for  an $E_V$-measurable, almost everywhere finite complex-valued function $\psi(\y),$ $\y\in\R^1,$ the operator $\psi(V)$ is defined as
\begin{equation}\label{psi(V)}
    \psi(V)=\sum_{\ell\in\Z_+}\psi(\ell)\Pb_\ell,
\end{equation}
with the natural domain
\begin{equation*}
    \Dc(\psi(V))=\{f=\sum f_\a e_{\a}\in\Ac_\l(\B^{\nb}): \sum_\ell |\psi(\ell)|^2\sum_{|\a|=\ell}|f_\a|^2<\infty\}.
\end{equation*}
On this domain, the operator $\psi(V)$ is normal; if $\psi$ has real values at all integer  points, $\psi(V)$ is self-adjoint; it is bounded iff $\psi$ is a bounded function on $\Z_+.$

Similarly to the case $\nb=1$ considered in \cite{RVSpectralDisk}, the operator $\psi(V)$ is compact iff  $\psi(\ell)\to 0$ as $\ell\to\infty.$ The explanation of this property is just a little bit more intricate than for $\nb=1$, since the projections $\Pb_\ell$ are not one-dimensional any more and their rank is not uniformly bounded.  Namely, for any given $\e>0$, let $\D(\e)$ be the set of those $\ell\in \Z_+$ where $|\psi(\ell)|>\e.$ This set is finite, moreover the spectral projection $E_V(\D(\e))$ is finite-dimensional (here we use the fact that all eigenspaces of $V$ are finite-dimensional). Thus the operator $\psi_\e(V)=\sum_{\ell\in \D(\e)}\psi(\ell)\Pb_\ell$ has finite rank, in particular, it is compact. On the other hand, for $\psi_\e'(\y)=\psi(\y)-\psi_\e(\y)$, we have $\|\psi_\e'(V)\|\le\e,$ since $|\psi(\y)|\le \e$ for $\y\notin \D(\e).$ So, $\psi(V)$ is $\e$-approximated in the operator norm by a compact operator and therefore it is compact itself.  The converse implication is also justified rather easily.

For different functions $\psi_1,\psi_2$ the operators $\psi_1(V), \psi_2(V)$ commute: in the strict sense, for bounded $\psi_1,\psi_2,$ or in the resolvent sense for more general functions, as in \cite{RVSpectralDisk}.

Being applied to the function $\psi=\pmb{1}$, $\pmb{1}(\ell)=1$ for all $\ell\in \Z_+,$ \eqref{psi(V)} gives $\pmb{1}(V)|_{\Ac_{\l}}=\mathbf{Id}_{\Ac_{\l}(\B\nb)},$ the identity operator in $\Ac_\l(\B^\nb).$ On the other hand, if $f\in \Hc_\l(\B^\nb)$ and $f$ is orthogonal to all functions in $\Ac_\l(\B^\nb)$, we have $\pmb{1}(V)f=0$ since each $\Pb_\ell f$ equals zero. Therefore, $\pmb{1}(V)$ is nothing else but the Bergman projection $\mathbb{P}\equiv\mathbb{P}_{\l}:\Hc_{\l}(\B^\nb)\to  \Ac_\l(\B^\nb).$

Let $a(r)=a(|z|)$ be a radial symbol on the ball $\B^\nb.$ Suppose that $a(r)\in L_\infty(0,1).$ Following \cite{RVSpectralDisk}, we establish a relation between radial Toeplitz operators in $\Ac_\l(\B^\nb)$ and functions of the rotation operator $V.$

The Toeplitz operator $\Tb_a$ in $\Ac_\l(\B^\nb)$ with symbol $a\in L_\infty$  acts as
\begin{equation}\label{Project}
    \Tb_a f =\mathbb{P} a f.
\end{equation}
In particular, for $f=e_{\a},$ \eqref{Project} gives
\begin{equation}\label{ProjectBas}
    (\Tb_a e_{\a})(z)=(\mathbb{P} a  e_{\a})(z)=\sum_{\ell}\Pb_{\ell}(a(|\z|)e_{\a}(\z))(z).
\end{equation}
Now we recall that the multiplication by a radial symbol retains the orthogonality of the basic functions $e_{\a}.$
Therefore, the operator  $\Tb_a$ is diagonal in the basis $\{e_{\a}\}.$ As well known (see, e.g., \cite{GKV}), in each eigenspace $\Ls_\ell,$ the elements on the diagonal of $\Tb_a$ depend only on $\ell,$ in particular, they are equal. It follows that the operator  $\Tb_a$ has block structure with respect to the eigenspaces $\Ls_\ell,$ with each block of dimension $\db_\ell$ having the form $\g_{\ell,a}\mathbf{Id}_{\ell},$ this latter symbol denoting the identity operator in $\Ls_{\ell}.$ This is exactly the structure we described above for  functions of the operator $V.$ Namely, with the radial symbol $a(r),$
we associate the function $\psi_a$
\begin{equation}\label{spectral}
    \psi_a(\ell)=\g_{\ell,a}, \ell\in\Z_+.
\end{equation}
Thus we obtain the description of a class of radial Toeplitz operators.

\begin{theorem}\label{ThmStructure} For a radial symbol $a,$ the Toeplitz operator $\Tb_a$ in $\Ac_\l(\B^\nb)$ is a function of the operator $V$,
\begin{equation*}
    \Tb_a=\psi_a(V).
\end{equation*}
The operator $\Tb_a$ is bounded iff the function $\psi_a$ is bounded on $\Z_+$; $\Tb_a$ is compact iff $\psi_a(\ell)\to 0$ as $\ell\to\infty.$
\end{theorem}

It follows from the expression for $\g_{\ell,a}$ that for a bounded radial symbol $a(r)$, the function  $\psi_a$ is bounded and therefore the operator $\Tb_a$ is bounded. The relation $a\Leftrightarrow \psi_a$ can be extended to a wide class of unbounded symbols $a.$ Since the expression for $\psi_a$ for the ball $\B^\nb$ is formally the same as for the one-dimensional complex disk $\B^1$ (up to a re-numeration), the results in \cite{RVSpectralDisk} carry over to the multidimensional case without changes.  We describe them here briefly.

On the one hand, bounded Toeplitz operators can be defined for symbols considerably more general than the bounded ones. First, as in    \cite{RVSpectralDisk}, quite a wide class of symbols with compact support is admissible. Any $a\in L_1(0,1)$ that vanish in some neighborhood of the endpoint $1$ generates a bounded function $\psi_a.$
Moreover, the same is correct for $a$ being a distribution with compact support in $[0,1).$ Although the reasoning in \cite{RVSpectralDisk} cannot be carried over directly to the multi-dimensional case (some complications in dealing with the point $0$ arise), the corresponding reasoning in \cite{RVBergDisk} works, namely by means of considering the symbol $a$ as a distribution in $\Es'(\B^{\nb}),$ invariant with respect to rotations of the ball. A class of symbols with support touching the point $1$ can be considered as well, including the ones unbounded in any neighborhood of $1$ but fast oscillating as $r\to 1$.

On the other hand, it is important to keep in mind that not for all bounded functions  $\psi(\ell)$ the operator $\psi(V),$ even a bounded one, is a radial Toeplitz operator with a bounded radial symbol $a.$ The necessary and sufficient condition for this can be derived from \cite{GKV}. Namely, the sequence $\psi(\ell)$ must be slowly oscillating, see \cite{RVSpectralDisk}.

The condition for the compactness of  Toeplitz operators,  expressed in terms of the function $\psi_a$, is given above: it requires that $\psi_a(\ell)\to 0$ as $\ell\to\infty.$ In terms of the symbol $a,$ this condition involves a certain (in some mean sense) vanishing of $a(r)$ as $r\to 1.$

\section{Toeplitz operators with quasi-radial symbols} \label{se:quasi-radial}

Now we pass to the algebras of Toeplitz operators with symbols of more complicated structure, namely the ones being radial in separate groups of variables. Each class of symbols is invariant with respect to a certain isometry group of the ball $\B^\nb$ described by a partition $\kn,$ as in Section \ref{isometry}. The classes considered in two previous sections serve as extreme cases.

Further on, in this Section,  the   partition $\kn=(k_1,...,k_m)$ is assumed fixed. A measurable function $a=a(z)$, $z \in \mathbb{B}^\nb$, will be called \emph{$\kn$-quasi-radial} if it depends only on $r_{(1)}$, ..., $r_{(m)}$, where
\begin{equation} \label{eq:part_rad}
 r_{(j)} = |z_{(j)}| = \sqrt{|z_{j,1}|^2 + \ldots + |z_{j,k_j}|^2}, \qquad j=1,\ldots, m.
\end{equation}

With this definition, such symbols can vary from \emph{separately radial}, $a = a(|z_1|, \ldots, |z_\nb|)$, if $\kb=(1,...,1)$, to \emph{radial}, $a = a(|z|)$, if $\kn=(\nb);$  these extreme cases have been considered in the previous sections.

Recall, see \cite[Lemma 3.1]{Vas_q-rad}, that the Toeplitz operator $\Tb_a$ in $\Ac_\l(\B^{\nb})$, with $\kn$-quasi-radial symbol $a=a(r_{(1)},\ldots,r_{(m)})$, acts on the basis elements $e_{\alpha}(z)$ as
\begin{equation} \label{eq:on_e_alpha}
 \Tb_a\, e_{\alpha} = \gamma_{a,\kn,\lambda}(\alpha)\,e_{\alpha}, \ \ \ \ \alpha \in \mathbb{Z}^\nb_+,
\end{equation}
where
\begin{eqnarray} \label{eq:gamma_quasi-ragial}
 && \gamma_{a,\kn,\lambda}(\alpha) =  \gamma_{a,\kn,\lambda}(|\alpha_{(1)}|,...,|\alpha_{(m)}|) \\ \nonumber
&& = \ \frac{2^m\, \Gamma(\nb+|\alpha| + \lambda + 1)}{\Gamma(\lambda+1)  \prod_{j=1}^m (k_j-1+|\alpha_{(j)}|)!}
 \int_{\tau(\mathbb{B}^m)} a(r_{(1)},...,r_{(m)}) (1-|r|^2)^{\lambda}  \prod^m_{j=1} r_{\!(j)}^{2|\alpha_{(j)}|+2k_j-1} dr_{\!(j)} \\ \nonumber
&& = \ \frac{ \Gamma(\nb+|\alpha| + \lambda + 1)}{\Gamma(\lambda+1)  \prod_{j=1}^m (k_j-1+|\alpha_{(j)}|)!}
 \int_{\pmb{\Delta}_{m}} a(\sqrt{\r_1},...,\sqrt{\r_m}) (1-(\r_{1}+\ldots + \r_m))^{\lambda}  \\ \nonumber
 && \times \ \prod^m_{j=1} \r_j^{|\alpha_{(j)}|+k_j-1} d\r_j,
\end{eqnarray}
and $\pmb{\Delta}_{m} = \{\r \in \mathbb{R}_+^m : \r_1+\ldots + \r_m < 1\}$ denotes the standard simplex, $\r_{j}=r_{(j)}^2$.

We will denote by $\mathcal{T}(\kn)$ the $C^*$-algebra, generated by Toeplitz operators with $\kn$-quasi-radial symbols.

Given a multi-index $\sbb = (s_1,\ldots,s_m) \in \mathbb{Z}_+^m$, we define (see e.g. \cite[Formula (3.5)]{BV_2013}) the finite dimensional subspace $H_\sbb$ in $\mathcal{A}_\l(\B^{\nb})$ as
\begin{equation*}
 H_\sbb = \mathrm{span}\,\{e_{\alpha}: \ |\alpha_{(j)}| = s_j,\quad j=1,\ldots,m\}.
\end{equation*}
This subspace consists of functions which are, for each $j=1,\dots,m,$   homogeneous polynomials of order $s_j$ in the group of variables $z_{(j)}$.
Then, given $\ell \in \mathbb{Z}_+$, for any fixed $j,$ we define also (see e.g. \cite[Formula (3.10)]{BV_2013}) the \emph{infinite dimensional} subspace $H^{(j)}_\ell$ in $\mathcal{A}_\l(\B^{\nb})$ as
\begin{equation*}
 H^{(j)}_\ell = \overline{ \mathrm{span}}\,\{e_{\alpha}: \ |\alpha_{(j)}| =\ell \}.
 \end{equation*}
 The latter space is the closure in $\Ac_\l(\B^{\nb})$ of the space of polynomials in $z$ which are order $\ell$ polynomials with respect to  the group of variables $z_{(j)}.$ With this definition,
 $H_\sbb = \bigcap_{j=1,\ldots,m} H^{(j)}_{s_j}$.
We introduce  the corresponding orthogonal projections
\begin{equation*}
 \Pb_\sbb : \mathcal{A}_\l(\B^{\nb})\longrightarrow H_\sbb \qquad \mathrm{and} \qquad \Pb^{(j)}_{\ell} : \mathcal{A}_\l(\B^{\nb}) \longrightarrow H^{(j)}_\ell, \quad \mathrm{with} \quad \Pb_\sbb = \prod_{j=1,\ldots,m} \Pb^{(j)}_{s_j}.
\end{equation*}

In this way, the Bergman space $\mathcal{A}_\l(\B^{\nb})$ splits into the orthogonal sum in two ways:
\begin{equation} \label{eq:orth_sum}
 \mathcal{A}_\l(\B^{\nb}) = \bigoplus_{\sbb \in \mathbb{Z}_+^m} H_\sbb
 \qquad \mathrm{and} \qquad \mathcal{A}_\l(\B^{\nb}) = \bigoplus_{\ell \in \mathbb{Z}_+} H^{(j)}_\ell.
\end{equation}
 We see that in the second case, each summand consists of functions having the fixed homogeneity order in variables in the group $z_{(j)}$ and analytic in the remaining variables.

\smallskip
Further, \eqref{eq:on_e_alpha} implies a decomposition of the Toeplitz operator:
\begin{equation} \label{eq:Toeplitz_quasi-radial}
 \Tb_a = \sum_{\sbb \in \mathbb{Z}_+^m} \gamma_{a,\kn,\lambda}(\sbb) \Pb_\sbb,
\end{equation}
understood in the sense of the strong convergence.

For each $j = 1,\ldots, m$, we introduce the rotation operator along the group of variables $z_{(j)}$:
\begin{equation} \label{eq:V_(j)}
 V_{(j)} = z_{j,1} \frac{\partial }{\partial z_{j,1}} + \ldots +
 z_{j,k_j} \frac{\partial }{\partial z_{j,k_j}} = \frac{1}{\imath} \frac{\partial }{\partial \theta_{j,1}} + \ldots + \frac{1}{\imath} \frac{\partial }{\partial \theta_{j,k_j}}, \quad \mathrm{where} \ \ z_{j,\ell} = |z_{j,\ell}|e^{\imath\theta_{j,\ell}},
\end{equation}
which is self-adjoint being defined on
\begin{equation*}
 \mathcal{D}(V_{(j)}) = \left\{f\in\Ac_\l(\B^{\nb}), f = \sum_{\alpha \in \mathbb{Z}_+^\nb} f_{\alpha}e_{\alpha} : \ \sum_{\alpha \in \mathbb{Z}_+^\nb}
 |\alpha_{(j)}|^2 |f_{\alpha}|^2 < \infty \right\}.
\end{equation*}

Note that on this domain, the operator $V_{(j)}$ admits the representation
\begin{equation} \label{eq:V_j}
 V_{(j)} = \sum_{\ell \in \mathbb{Z}_+} \ell \Pb^{(j)}_{\ell},
\end{equation}
understood in the sense of the strong convergence.

For each $j$ the operator $V_{(j)}$ acts only upon variables $z_l$ entering in  $z_{(j)}$ and is the identity operator in all other variables. Therefore, these operators commute.
For each fixed $j,$ the operator $V_{(j)}$ has spectral structure similar to the one we described in Sect. \ref{Sect.radial} for a ball, with the only difference that this ball  has dimension $k_j$ instead of $\nb.$ Therefore,
  the spectrum of $V_{(j)}$ consists of isolated points $\ell\in\Z_+$. Namely, the spectral subspace   of $V_{(j)}$, corresponding to the eigenvalue $\ell$ is spanned by the eigenfunctions  of the form
  \begin{equation}\label{calcul.quasi}
  f_\ell^{(j)}(z)=z_{(j)}^{\a_{(j)}}h(\widetilde{z_{(j)}})
\end{equation}
with some multiindex $\a_{(j)}\in\Z_+^{k_j},$
$|\a_{(j)}|=\ell$  and an analytic function $h$ of the complementing variables $\widetilde{z_{(j)}}.$
The calculation, similar to the one in \eqref{calculation2} shows that the function $f_\ell^{(j)}(z)$ of the form \eqref{calcul.quasi} belongs to the Bergman space $\Ac_\l(\B^\nb)$ if and only if the function $h(\widetilde{z_{(j)}})$ of the complementing variables  belongs to the Bergman space $\Ac_{\l+\ell+k_j}(\B^{\nb-k_j})$ in the ball of a smaller dimension $\nb-k_j$, with equivalence of norms.

    Finite linear combinations of elements in different subspaces $H^{(j)}_\ell$ are dense in $\Ac_\l(\B^{\nb})$, therefore there are no other points of spectrum of $V_{(j)}$. Each point $\ell\in\Z_+$ is thus an eigenvalue of $V_{(j)}$ of \emph{infinite multiplicity},  and $\Pb^{(j)}_{\ell}$ is the spectral projection of $V_{(j)}$ corresponding to the eigenvalue $\ell\in\Z_+$.

The Spectral Theorem representation of the operator $V_{(j)}$ is therefore
\begin{equation*}
 V_{(j)} = \int_{\mathbb{R}} \y_j d\Eb_j(\y_j),
\end{equation*}
where the resolution of identity $\Eb_j(\eta_j) \equiv E_j((-\infty,\eta_j])$ is given by
\begin{equation*}
 \Eb_j(\eta_j) = \sum_{\ell \leq \eta} \Pb^{(j)}_{\ell}.
\end{equation*}

\medskip
Having  done this for all $j=1,2,...,m$, we construct the joint spectral measure on $\R^m$  using $E_1$, $E_2$, ..., $E_m$, as described in Sect. 2. The measure $E$  is supported on the integer lattice $\Z^m_+$, while the set $\R^m \setminus \Z^m_+$ has zero $E$-measure. Note that, for each $j=1,2,...,m$, we have also
\begin{equation*}
 V_{(j)} = \int_{\mathbb{R}^m} \eta_j d\Eb(\eta_1,...,\eta_m),
\end{equation*}
where $\Eb(\by)$ is the resolution of identity corresponding to the measure $E$,
\begin{equation*}
  \Eb(\by)=\prod_{j}\Eb_j(\y_j).
\end{equation*}

The Functional Calculus implies now
\begin{proposition}
 Given an $E$-measurable almost everywhere finite function $\psi$ on $\mathbb{R}^m$, the operator
 \begin{equation*}
  \psi(\Vb)=\psi(V_{(1)},\ldots, V_{(m)}) := \int_{\mathbb{R}^m} \psi(\eta_1,\ldots, \eta_m)\,d\Eb(\eta_1,\ldots, \eta_m)
 \end{equation*}
is well defined and normal on its domain
\begin{equation*}
 \mathcal{D}_{\psi} = \left\{f \in \Ac_{\l}(\B^{\nb}) : \int_{\mathbb{R}^m} |\psi(\by)|^2\, d\langle \Eb(\by)f,f\rangle < \infty \right\}, \qquad \by=(\y_1,\dots,\y_m).
\end{equation*}
Moreover the operator $\psi(\Vb)$ is bounded, and thus defined on the whole $\Ac_{\l}(\B^{\nb})$, if and only if the function   $\psi$ is $E$-essentially bounded.
\end{proposition}

Representation \eqref{eq:calculus} implies that
\begin{equation*}
 \psi(\mathbf{V}) \equiv \psi(V_{(1)},\ldots, V_{(m)}) = \int_{\sigma(\mathbf{V})} \psi(\by)\,d\Eb(\by) = \sum_{\sbb \in \Z^m_+} \psi(\sbb)\Pb_\sbb.
\end{equation*}

\begin{corollary} \label{co:R(k)}
 The set of all operators $\{\psi(\mathbf{V})\}$, defined by (the classes of equivalence) of  $E$-measurable essentially bounded functions  $\psi$ with $\pmb{\pmb{\psi}}= \{\psi(\sbb)\}_{\sbb \in \mathbb{Z}^m_+} \in \ell_{\infty}(\mathbb{Z}^m_+)$ constitutes the algebra $\mathcal{R}(\kn)$ of bounded mutually commuting operators in $\mathcal{A}_{\lambda}(\mathbb{B}^{\nb})$, and the mapping
 \begin{equation*}
  Op_{E}:\,\boldsymbol\psi \  \longmapsto \psi(\mathbf{V}) = \sum_{\sbb \in \Z^m_+} \psi(\sbb)\Pb_\sbb.
 \end{equation*}
defines the isomorphism
\begin{equation*}
 Op_{E}:\, \ell_{\infty}(\mathbb{Z}^m_+) \ \longrightarrow \ \mathcal{R}(\kn)\,\subset\Bc(\Ac_\l(\B^{\nb}))
\end{equation*}
of commutative $C^*$-algebras.
\end{corollary}

We return now to Toeplitz operators with $\kn$-quasi-radial symbols. By \eqref{eq:Toeplitz_quasi-radial}, a Toeplitz operator $\Tb_a$, with $\kn$-quasi-radial symbol $a= a(r_{(1)},r_{(2)},\ldots,r_{(m)})$, admits the representation
\begin{equation*}
 \Tb_a = \sum_{\sbb \in \mathbb{Z}_+^m} \gamma_{a,\kn,\lambda}(\sbb) \Pb_\sbb,
\end{equation*}
where the coefficients $\gamma_{a,\kn,\lambda}(\sbb)$ are given by \eqref{eq:gamma_quasi-ragial}. Thus, by Corollary \ref{co:R(k)},
\begin{equation} \label{eq:T_a-spectral}
 \Tb_a = \psi_a(V_{(1)},V_{(2)},\ldots,V_{(m)}),
\end{equation}
where $\psi_a$ belongs to the equivalence class of essentially bounded $E$-measurable functions on $\R^m$, defined by the sequence $\pmb{\pmb{\psi}}_a= \{\psi_a(\sbb)\}_{\sbb \in \mathbb{Z}^m_+}$, with $\psi_a(\sbb) = \gamma_{a,\kn,\lambda}(\sbb)$.

\medskip
It seems to be useful to enrich the situation by adding some important properties and characteristics of the objects under study in terms  of the representation theory. \\ Our considerations are based upon  \cite{Quiroga_21}, where all definitions, details, and proofs can be found.

Recall, that a unitary representation $\pi_{\lambda}: \mathbf{U}(\nb) \rightarrow \mathcal{B}(\mathcal{A}_{\lambda}(\B^\nb))$ of the group $\mathbf{U}(\nb)$ of all $n \times n$ unitary matrices on the weighted Bergman space $\mathcal{A}_{\lambda}(\B^\nb)$ is given by
\begin{equation*}
 \pi_{\lambda}(U)(f) = f \circ U^{-1}, \, U\in \Ub(\nb).
\end{equation*}
For a closed subgroup $\Hb$ of $\mathbf{U}(\nb)$, we denote by $\pi_{\lambda}|_\Hb$ the restriction of the representation $\pi_{\lambda}$ to the group $\Hb$. Then, an operator $T \in \mathcal{B}(\mathcal{A}_{\lambda}(\B^\nb))$ intertwines the restriction $\pi_{\lambda}|_\Hb$, equivalently $T$ is $\Hb$-equivariant, if
\begin{equation*}
 T\pi_{\lambda}(U) = \pi_{\lambda}(U)T, \quad \mathrm{for \ all} \quad U \in \Hb.
\end{equation*}
We will denote by $\mathrm{End}_\Hb (\mathcal{A}_{\lambda}(\B^\nb))$ the algebra of all such intertwining, or $\Hb$-equivariant, operators. It is well known that $\mathrm{End}_\Hb (\mathcal{A}_{\lambda}(\B^\nb))$ is a von Neumann algebra.

For our  partition $\kn=(k_1,...,k_m)$, let
\begin{equation*}
 \Hb = \mathbf{U}(\kn) := \mathbf{U}(k_1) \times \ldots
 \times \mathbf{U}(k_m),
\end{equation*}
realized as a block diagonal subgroup of $\mathbf{U}(\nb)$.

Note that \cite[Proposition 4.5]{Quiroga_21} demonstrates the importance of the first orthogonal sum decomposition in \eqref{eq:orth_sum}. We give an equivalent formulation of this proposition, in which the summands are defined in terms of basis elements of $\mathcal{A}_{\lambda}(\B^\nb))$, while the original formulation defines them in terms of certain polynomials in $z$.

\begin{proposition}[{\cite[Proposition 4.5]{Quiroga_21}}]
 Given a partition $\kn=(k_1,...,k_m)$, the isotypic decomposition of the restriction $\pi_{\lambda}|_{\mathbf{U}(\kn)}$ is
 \begin{equation} \label{eq:isotypic}
 \mathcal{A}_{\lambda}(\B^\nb) = \bigoplus_{\sbb \in \mathbb{Z}_+^m} H_\sbb.
 \end{equation}
Furthermore, this isotypic decomposition is multiplicity-free.
\end{proposition}
This means that the restriction $\pi_{\lambda}|_{\mathbf{U}(\kn)}$ onto each summand $\Hb_\sbb$ is irreducible, and each two summands in the above orthogonal sum decomposition are non-isomorphic irreducible $\mathbf{U}(\kn)$-modules. Then by the Schur Lemma,  each $\mathbf{U}(\kn)$-equivariant operator preserves isotypic components in \eqref{eq:isotypic}, acting on each $\Hb_\sbb$ as the multiplication by a scalar operator. This means that each operator
$T \in \mathrm{End}_{\mathbf{U}(\kn)} (\mathcal{A}_{\lambda}(\B^\nb))$ must have the form
\begin{equation*}
 T = \sum_{\sbb \in \mathbb{Z}_+^m} \gamma(\sbb) \Pb_\sbb,
\end{equation*}
where, recall, $\Pb_\sbb: \mathcal{A}_{\lambda}(\B^\nb)) \rightarrow \Hb_\sbb$ is the orthogonal projection.

An equivalent reformulation of \cite[Proposition 4.9]{Quiroga_21} says now
\begin{proposition}[{\cite[Proposition 4.9]{Quiroga_21}}] \label{prop:End}
Given a partition $\kn=(k_1,...,k_m)$, the mapping
\begin{equation*}
\ell_{\infty}(\mathbb{Z}_+^m) \ \longrightarrow \ \mathrm{End}_{\mathbf{U}(\kn)} (\mathcal{A}_{\lambda}(\B^\nb))
\end{equation*}
defined by
\begin{equation*}
 \{\gamma(\sbb)\}_{\sbb \in \mathbb{Z}_+^m} \ \longmapsto \ \sum_{\sbb \in \mathbb{Z}_+^m} \gamma(\sbb) \Pb_\sbb
\end{equation*}
is the isomorphism of the $C^*$-algebras which are  von Neumann algebras at the same time.
\end{proposition}

The results of Corollary \ref{co:R(k)} and Proposition \ref{prop:End} lead now to the following statement.

\begin{proposition}
 The $C^*$-algebra $\mathcal{R}(\kn) = Op_{E}( \ell_{\infty}(\mathbb{Z}^m_+))$ coincides with the von Neumann algebra $\mathrm{End}_{\mathbf{U}(\kn)} (\mathcal{A}_{\lambda}(\B^\nb))$ of all bounded $\mathbf{U}(\kn)$-equivariant operators. Each operator $T \in \mathrm{End}_{\mathbf{U}(\kn)} (\mathcal{A}_{\lambda}(\B^\nb))$ is thus a certain function of the commuting operators  $V_1,\ V_2,\ \ldots, \ V_m$,
 \begin{equation*}
  T = \psi(V_{(1)},V_{(2)},\ldots, V_{(m)}),
 \end{equation*}
where $\psi$ belongs to the equivalence class of essentially bounded $E$-measurable functions, defined by the sequence $\pmb{\boldsymbol\psi}= \{\psi(\sbb)\}_ {\sbb \in \mathbb{Z}^m_+}$, with
\begin{equation*}
 \psi(\sbb) = \frac{\mathrm{tr}(T|_{H_{\sbb}})}{\mathrm{dim}\,H_{\sbb}}.
\end{equation*}
\end{proposition}

\medskip
Now instead of a fixed partition $\kn=(k_1,k_2,\ldots,k_m)$, we consider the finite set $\K$ of all possible partitions of $\nb$ into positive integers. We introduce  the partial order on $\K$ as
\begin{equation*}
 \kn'=(k'_1,k'_2,\ldots,k'_{m'}) \  \preccurlyeq \ \kn''=(k''_1,k''_2,\ldots,k''_{m''})
\end{equation*}
if and only if for each $j =1,2,\ldots,m'$, $k'_j$ is a sum, after a re-numeration, of certain consecutive elements of the partition $\kn''$,
\begin{equation*}
 k'_j = k''_{p_j} + k''_{p_j+1} + \ldots + k''_{p_j+q_j-1},
\end{equation*}
or if and only if for each $j =1,2,\ldots,m$, after a reordering,
\begin{equation*}
 z'_{(j)} = (z''_{p_j}, z''_{p_j+1}, \ldots, z''_{p_j+q_j-1}),
\end{equation*}
equivalently if and only if $\mathbf{U}(\kn'') \subseteq \mathbf{U}(\kn')$.

With respect to  this  partial order, $\kn = (\nb)$ and $\kn = (1,1,\ldots,1)$ are the minimal and the maximal elements in~$\K$, respectively. Furthermore, this order implies the order by inclusion of the von Neumann algebras $\mathrm{End}_{\mathbf{U}(\kn)} (\mathcal{A}_{\lambda}(\B^\nb))$ of the sets of bounded $\mathbf{U}(\kn)$-equivariant operators:
\begin{equation*}
 \mathrm{End}_{\mathbf{U}(\kn')} (\mathcal{A}_{\lambda}(\B^\nb)) \ \subseteq \ \mathrm{End}_{\mathbf{U}(\kn'')} (\mathcal{A}_{\lambda}(\B^\nb)) \quad \mathrm{if \ and \ only \ if} \quad \kn' \  \preccurlyeq \ \kn''.
\end{equation*}
Note that the  subalgebras $\mathcal{T}(\kn)$ of $\mathrm{End}_{\mathbf{U}(\kn)} (\mathcal{A}_{\lambda}(\B^\nb))$, generated by Toeplitz operators with bounded symbols, maintain the same order,
\begin{equation*}
 \mathcal{T}(\kn') \ \subseteq \ \mathcal{T}(\kn'') \quad \mathrm{if \ and \ only \ if} \quad \kn' \  \preccurlyeq \ \kn''.
\end{equation*}

\section{Degenerate elliptic commutative $C^*$-algebras}

In the preceding sections we have demonstrated that the size of the set of commuting operators involved in the spectral representation of Toeplitz operators with a certain class of symbols depends on the size of this set of symbols: the larger is the algebra of symbols, the more operators are needed for the spectral representation. In this last section we discuss two examples demonstrating to what extent this rule applies to classes of symbols in a more general setting.

 The symbols to be consider are characterized by their (essential) dependence  only on  a part of variables $z_{(j)}$ (or only on a non-complete set of their   combinations), therefore, we call them \emph{degenerate elliptic}. Such symbols and the corresponding Toeplitz operators were studied, in particular, in \cite{BHaV,BV_2017}, see also \cite[Section 6]{Quiroga_Sanchez_JFA}, where they  are present under the name $\b$-quasi-elliptic symbols.

As before,  we rearrange (via an appropriate biholomorphism of $\mathbb{B}^{\nb}$) the coordinates of $z=(z_1,...,z_\nb) \in \mathbb{B}^{\nb}$ in the most convenient order for the further description. Thus, without loss of generality, we will use the following setup.

Let, as Section \ref{se:quasi-radial}, $\kn=(k_1,k_2,...,k_m)$ be a partition of $\nb$, with correspondingly $z=(z_{(1)},\dots,z_{(m-1)},z_{(m)})$. We denote by $\kn'$ the partition $\kn'=(k_1,k_2,...,k_{m'}),$ $m'=m-1$ of the integer $\nb'=\nb-k_{m}.$ The symbols under consideration depend on the first $m'$ radii $r_{(j)}=z_{(j)},$ $j=1,\dots,m-1,$ where $r'=(r_{(1)},\dots,r_{(m-1)}).$  The dependence on  $r''\equiv r_{(m)}=|z_{(m)}|$ is realized, in our first example, via the presence of the weight,
  \begin{equation}\label{degen.weight}
 a_w(z)=a_w(r',r'')=a\left(\frac{r_{(1)}}{\sqrt{1 - |r''|^2}}, \ldots, \frac{r_{(m')}}{\sqrt{1 - |r''|^2}}\right).
\end{equation}
Here, although all variables are present,  the symbol is, actually,  a function of $m'$ variables via the superposition.

In the second example, the  symbols are independent of $z_{(m)}:$
\begin{equation}\label{degen}
  a(z)=a(r_{(1)},\dots,r_{(m-1)} )=a(r').
\end{equation}
Of course, both classes of symbols belong to the set of $\kn$-quasi-radial ones, so the general spectral representation involving $m$ operators $V_{(j)}$ is valid. Since each of these classes is considerably smaller than the algebra of $\kn$-quasi-radial symbols, one might expect that the spectral representation  uses fewer operators.



\subsection{Weighted $\kn'$-quasi-radial symbols}

We consider the algebra $\SF_w(\kn')$ of all bounded symbols $a_w(z)$ in \eqref{degen.weight} and call them \emph{weighted $\kn'$-quasi-radial symbols}.

Since each weighted $\kn'$-quasi-radial symbol is quasi-radial for our $\kn =(\kn',k_m)$, the corresponding Toeplitz operator is diagonal,
\begin{equation*}
 \Tb_{a_w} e_\a(z)=\g_{a_w,\kn',\l} e_\a(z),
\end{equation*}
where, by \eqref{eq:gamma_quasi-ragial}, for a multi-index $\a=(\a_{(1)},\dots,\a_{(m-1)},\a'')\in \Z^\nb_+$,
\begin{eqnarray*}
&& \gamma_{a_w,\kn',\lambda}(\alpha) = \frac{2^m\, \Gamma(\nb+|\alpha| + \lambda + 1)}{\Gamma(\lambda+1)  \prod_{j=1}^{m'} (k_j-1+|\alpha_{(j)}|)! (k_m-1+|\alpha''|)!} \\
 && \times \
 \int_{\tau(\mathbb{B}^m)} a_w(r_{(1)},...,r_{(m')}, |r''|) (1-|r|^2)^{\lambda}  \prod^{m-1}_{j=1} {r_{(j)}}^{2|\alpha_{(j)}|+2k_j-1} dr_{(j)} {r''}^{2|\alpha''+2k_m-1} dr'' \\
 && = \ \frac{2^m\, \Gamma(\nb+|\alpha| + \lambda + 1)}{\Gamma(\lambda+1)  \prod_{j=1}^{m-1} (k_j-1+|\alpha_{(j)}|)! (k_m-1+|\alpha'')!}  \int_0^1  {r''}^{2|\alpha''|+2k_m-1} dr'' \,\cdot \, I,
\end{eqnarray*}
with
\begin{equation*}
 I =\int_{|r'|^2<1-|r''|^2} a_w(r_{(1)},...,r_{(m-1)}, |r''|) (1-|r'|^2-|r''|^2)^\l \prod^{m-1}_{j=1} {r_{(j)}}^{2|\alpha_{(j)}|+2k_j-1} dr_{(j)}.
\end{equation*}
Making the change of variables $r_{(j)} = \sqrt{1-|r''|^2}u_j$, $|u| =\sqrt{u_1^2+...+u_{m-1}^2}$, $\u_j=u_j^2$ we have
\begin{eqnarray*}
 && I = \int_{\tau(\mathbb{B}^{m-1})} a(u_{(1)},...,u_{(m-1)})(1-|r''|^2)^{\l+ |\a'|+ \nb'}(1-|u|^2)^\l \prod^{m-1}_{j=1} {u_j}^{2|\alpha_{(j)}|+2k_j-1} du_j \\
 && =2^{1-m}(1-|r''|^2)^{\l+ |\a'|+ \nb'}
 \int_{\pmb{\Delta}_{m'}} a(\sqrt{\u_1},...,\sqrt{\u_{m'}}) (1-\u_1-...-\u_{m-1})^\lambda
 \prod^{m'}_{j=1} {\u_j}^{|\alpha_{(j)}|+k_j-1} d\u_j.
\end{eqnarray*}
Thus, finally,
\begin{eqnarray*}
 && \gamma_{a_w,\kn',\lambda}(\alpha) = \frac{2^{-1}\, \Gamma(\nb+|\alpha| + \lambda + 1)}{\Gamma(\lambda+1)  \prod_{j=1}^{m-1} (k_j-1+|\alpha_{(j)}|)! (k_m-1+|\alpha''|)!} \\
 && \times \ \int_0^1 (1-|r''|^2)^{\l+ |\a'|+ \nb'} {r''}^{2|\alpha''+2k_m-1} dr'' \\
 && \times \ \int_{\pmb{\Delta}_{m''}} a(\sqrt{\u_1},...,\sqrt{\u_{m-1}}) (1-\u_1-...-\u_{m-1})^\lambda
 \prod^{m-1}_{j=1} {\u_j}^{|\alpha_{(j)}|+k_j-1} d\u_j \\
 && = \ \frac{\Gamma(\nb'+|\alpha'| + \lambda + 1)}{\Gamma(\lambda+1)  \prod_{j=1}^{m-1} (k_j-1+|\alpha_{(j)}|)!} \\
 && \times \ \int_{\pmb{\Delta}_{m-1}} a(\sqrt{\u_1},...,\sqrt{\u_{m-1}}) (1-\u_1-...-\u_{m-1})^\lambda
 \prod^{m-1}_{j=1} {\u_j}^{|\alpha_{(j)}|+k_j-1} d\u_j.
\end{eqnarray*}
The  most important result of these calculations is that the coefficient $\gamma_{a_w,\kn',\lambda}(\alpha)$ does not in fact depend on the last component $\a''=\a_{m}$ since all factors depending on $\a''$ cancel,
\begin{equation}\label{eq:gamma_w}
 \gamma_{a_w,\kn',\lambda}(\alpha) = \gamma_{a_w,\kn',\lambda}(|\alpha_{(1)}|,...,|\alpha_{(m-1)}|).
\end{equation}
Therefore, following Section \ref{se:quasi-radial}, for a given multi-index $\sbb = (s_1,\ldots,s_{m-1}) \in \mathbb{Z}_+^{m-1}$, we can define the, now,  the \emph{infinite dimensional} subspace $H_\sbb$ in $\mathcal{A}_\l(\B^{\nb})$  as
\begin{equation*}
 H_\sbb = \overline{\mathrm{span}}\,\{e_{\alpha}: \ |\alpha_{(j)}| = s_j,\quad j=1,\ldots,m-1 \}.
\end{equation*}
Let $\Pb_\sbb : \mathcal{A}_\l(\B^{\nb})\longrightarrow H_\sbb$ be the orthogonal projection.

In this way, given a weighted $\kn'$-quasi-radial symbol $a_w$, the  equality \eqref{eq:Toeplitz_quasi-radial} implies the representation
\begin{equation*}
 \Tb_{a_w} = \sum_{\sbb \in \mathbb{Z}_+^{m'}} \gamma_{a_w,\kn',\lambda}(\sbb) \Pb_\sbb,
\end{equation*}
understood in the sense of the strong convergence. The difference with the, formally similar, expression in Section \ref{se:quasi-radial} consists in the fact that now the projections $\Pb_\sbb$ are infinite-dimensional.

The rotation operators $V_{(j)}, \, j=1,\dots,m-1$ are defined in the same way as before, they commute, and we arrive at  the spectral representation using a smaller set of generating operators, just as we expected.
\begin{equation} \label{eq:spectral_w}
 \Tb_a = \psi_a(V_{(1)},V_{(2)},\ldots,V_{(m')}),
\end{equation}
where $V_{(j)}$ are given by \eqref{eq:V_(j)} and $\psi_a$ belongs to the equivalence class of essentially bounded $E$-measurable functions on $\R^{m'}$, defined by the sequence $\pmb{\boldsymbol\psi}_a= \{\psi_a(\sbb)\}_{\sbb \in \mathbb{Z}^{m'}_+}$, with $\psi_a(\sbb) = \gamma_{a_w,\kn',\lambda}(\sbb)$.

An \emph{essential} difference is that the reasoning in Sect.3 proving the compactness of the Toeplitz operator under the assumption that $\psi(\sbb)\to 0$ as $|\sbb|\to 0$ breaks down, again, due to the  infiniteness of the multiplicity of each point of the joint spectrum. Thu, such Toeplitz operators can never be compact.

\subsection{$\kn'$-quasi-radial symbols}
Using the same notations as in the Section 6.1, we consider now the symbols satisfying \eqref{degen}. We denote by $\SF(\kn')$ the algebra of such symbols.

For a multiindex $\a=(\a';\a_{(m)})=(\a_{(1)}, \dots,\a_{(m-1)};\a'' )$ we calculate, again, the action of the Toeplitz operator with $\kn'$-quasi-radial symbol upon the basic element $e_\a$. This calculation follows  the same way as for general quasi-radial symbols, however, some simplifications are present.

Toeplitz operators with symbols  $a\in\SF(\kn')$ are diagonal in the basis $\{e_\a(z\})=\{\e_{\a',\a''}(z',z'')\}$ and act on such  basis elements as
\begin{equation*}
 \Tb_a e_\a(z)=\g_{a,\kn',\l} e_\a(z),
\end{equation*}
where, by \eqref{eq:gamma_quasi-ragial},
\begin{eqnarray*}
 && \gamma_{a,\kn',\lambda}(\alpha) = \frac{2^{\nb'}\, \Gamma(\nb+|\alpha| + \lambda + 1)}{\Gamma(\lambda+1)  \prod_{j=1}^{m'} (k_j-1+|\alpha_{(j)}|)! (k_m-1+|\alpha''|)!} \\
 && \times
 \int_{\tau(\mathbb{B}^{\nb})} a(r_{\!(1)},...,r_{\!(m-1)}) (1-|r|^2)^{\lambda}  \prod^{m-1}_{j=1} {r_{(j)}}^{2|\alpha_{(j)}|+2k_j-1} dr_{\!(j)}  \, r''^{2|\alpha''|+2k_m-1} dr'' \\
 && = \ \frac{2^{\nb'}\, \Gamma(\nb+|\alpha| + \lambda + 1)}{\Gamma(\lambda+1)  \prod_{j=1}^{m'} (k_j-1+|\alpha_{(j)}|)! (k_m-1+|\alpha''|)!} \\
 && \times \ \int_{\tau(\mathbb{B}^{m-1})}a(r_{\!(1)},...,r_{\!(m-1)}) \prod^{m-1}_{j=1} {r_{\!(j)}}^{2|\alpha_{(j)}|+2k_j-1} dr_{\!(j)}\, \cdot \, I_{\a''}(|r'|),
\end{eqnarray*}
with
\begin{equation*}
I = \int_{|r''|^2<1-|r'|^2} (1-|r'|^2-|r''|^2)^\l {r''}^{2|\alpha''|+2k_m-1} dr''.
\end{equation*}
Making the change of variables $r''= \sqrt{1-|r'|^2}u$, we have
\begin{eqnarray*}
 && I = \int_0^1 (1-|r'|^2)^{\l+ |\a''|+ k_m}(1-|u|^2)^\l u^{2|\alpha''|+2k_m-1} du = \\
 && 2^{-1}(1-|r'|^2)^{\l+ |\a''|+ k_m}
 \int_0^1 (1-\u)^\lambda
 \u^{|\alpha''|+k_m-1} d\u=
  \\
 && 2^{-1}(1-|r'|^2)^{\l+ |\a''|+ k_m}\,
 \frac{\Gamma(\l+1)(|\alpha''|+k_m-1)!}{\Gamma(\l+|\a''|+k_m+1)}.
\end{eqnarray*}
Thus, finally,
\begin{eqnarray*}
 && \gamma_{a,\kn',\lambda}(\alpha) = \frac{2^{m'}\, \Gamma(\nb+|\alpha| + \lambda + 1)}{\Gamma(\l+|\a''|+\nb-\nb'+1)  \prod_{j=1}^{m'} (k_j-1+|\alpha_{(j)}|)!} \\
 && \times \ \int_{\tau(\mathbb{B}^{m-1})}a(r_{(1)},...,r_{(m-1)}) (1-|r'|^2)^{\l+ |\a''|+ \nb''} \prod^{m=1}_{j=1} {r_{(j)}}^{2|\alpha_{(j)}|+2k_j-1} dr_{(j)} \\
 && = \ \frac{\Gamma(\nb+|\alpha| + \lambda + 1)}{\Gamma(\l+|\a''|+\nb-\nb'+1)  \prod_{j=1}^{m-1} (k_j-1+|\alpha_{(j)}|)!} \\
 && \times \ \int_{\tau(\Delta_{m-1})}a(\sqrt{\r_{(1)}},...,\sqrt{\r_{(m-1)}}) (1-\r_{(1)}...-\r_{(m-1)})^{\l+ |\a''|+ \nb''} \prod^{m-1}_{j=1} {\r_{(j)}}^{|\alpha_{(j)}|+k_j-1} d\r_{(j)}.
\end{eqnarray*}
Observe that, in fact,
\begin{equation}\label{eq:gamma_m}
 \gamma_{a,\kn',\lambda}(\alpha) = \gamma_{a,\kn',\lambda}(|\alpha_{(1)}|,...,|\alpha_{(m-1)}|,|\alpha''|)
\end{equation}

Following Section \ref{se:quasi-radial}, for a given multi-index $\sbb = (s_1,\ldots,s_{m-1}, s_m) \in \mathbb{Z}_+^{m}$, we define the finite dimensional subspace $H_\sbb$ in $\mathcal{A}_\l(\B^{\nb})$ as
\begin{equation*}
 H_\sbb = \mathrm{span}\,\{e_{\alpha}: \ |\alpha_{(j)}| = s_j,\quad j=1,\ldots,m', \ \ |\a''| = s_{m} \},
\end{equation*}
and let $\Pb_\sbb : \mathcal{A}_\l(\B^{\nb})\longrightarrow H_\sbb$ be the orthogonal projection.

Further, given a $\kn'$-quasi-radial symbol $a$,  \eqref{eq:Toeplitz_quasi-radial} implies the representation
\begin{equation}
 \Tb_a = \sum_{\sbb \in \mathbb{Z}_+^{m}} \gamma_{a,\kn',\lambda}(\sbb) \Pb_\sbb,
\end{equation}
understood in the sense of the strong convergence.
Then Corollary \ref{co:R(k)} together with \eqref{eq:T_a-spectral} implies that
\begin{equation} \label{eq:spectral_m}
 \Tb_a = \psi_a(V_{(1)},V_{(2)},\ldots,V_{(m)}),
\end{equation}
where $V_{(j)}$ are given by \eqref{eq:V_(j)} and $\psi_a$ belongs to the equivalence class of essentially bounded $E$-measurable functions on $\R^m$, defined by the sequence $\boldsymbol\psi_a= \{\psi_a(\sbb)\}_{\sbb \in \mathbb{Z}^m_+}$, with $\psi_a(\sbb) = \gamma_{a,\kn',\lambda}(\sbb)$.

This calculation shows that for the spectral representation of Toeplitz operators with symbol in $\SF(\kn')$ we need the same system of $m$ operators $V_{(j)}$ as for the algebra of general quasi-radial symbols. This quantity cannot be reduced since the coefficients $\gamma_{a,\kn',\lambda}(\a)$ depend on all components of the multi-index $\a$, including $\a''$, unlike the case of weighted $\kn'$-quasi-radial symbols, where the coefficients are independent of $\a''$. This dependence is quite implicit, therefore it seems to be impossible to describe in a transparent way the set of those  sequences $\g(\a)$ which can serve as the function $\psi$ for this class of symbols. The question of the compactness of Toeplitz operators of this class seems to be rather     hard as well.
\subsection{Discussion: a more general setting}

To explain the  the difference of the results in two cases presented in the section, we embed them in a more general common setting. In both cases we started with a function $a$ depending on $\kn'= m-1$ radii $r_{(1)}$, ... ,$r'_{(m-1)}$. Formally, the  qualitative difference between these two cases is the following. In the first case, we include the dependence of $|z''|$ in the arguments of $a$. Here, both the resulting sequence $\gamma$ \eqref{eq:gamma_w} and the spectral representation of the corresponding Toeplitz operator \eqref{eq:spectral_w} depend only on the first group of $\kn'= m-1$ arguments, $|\alpha_{(1)}|$, ..., $|\alpha_{(m-1)}|$ and $V_{(1)}$, ..., $V_{(m-1)}$, respectively. In the second case, the symbol keeps depending only on $m-1$ arguments, however, both the resulting sequence $\gamma$ \eqref{eq:gamma_m} and the spectral representation of the corresponding Toeplitz operator \eqref{eq:spectral_m} depend on the whole  group of arguments,
$|\alpha_{(1)}|$, ..., $\alpha_{(m-1)}|$, $|\alpha"|$ and $V_{(1)}$, ..., $V_{(m)}$, respectively.

To understand the reason for such difference, it is instructive to combine the above two cases into one, more general, setup. We start with a partition $\kn=(\kn',\kn'')$ of  the number $\nb,$ such that the partition  $\kn'$ consists of $m'$ elements, $k'_1+ ... +k'_{m'}=\nb'$ and $\kn''$ consists of $m''$ elements, $k''_1+ ... +k''_{m''}=\nb''$, with $\nb'+ \nb'' = \nb$. Correspondingly, we divide the coordinates $z=(z_1,...,z_n) \in \mathbb{B}^{\nb}$ into groups, $z = (z',z'')$, where
\begin{eqnarray*}
 z' = \big{(}z'_{(1)}, \ldots, z'_{(m')}\big{)}, &\text{with}& z'_{(j)}
=\big{(} z'_{k'_1+\ldots +k'_{j-1}+1}, \ldots , z'_{k'_1+ \ldots +k'_j}\big{)}\in \mathbb{C}^{k'_j}\\
z'' = \big{(}z''_{(1)}, \ldots, z''_{(m'')}\big{)}, &\text{with}& z''_{(j)}
=\big{(} z''_{k''_1+\ldots +k''_{j-1}+1}, \ldots , z''_{k''_1+ \ldots +k''_j}\big{)}\in \mathbb{C}^{k''_j}.
\end{eqnarray*}
The following considerations will be based on the approach of \cite[Sections 2 and 3]{BHaV}, where more  details can be found.

Recall that the standard orthonormal basis of $\mathcal{A}_{\lambda}^2(\mathbb{B}^{\nb})$ is formed, see \eqref{basis}, by the weighted monomials
\begin{equation*}
e^{\lambda}_{\alpha}(z)= \sqrt{\frac{\Gamma(\nb+|\alpha|+\lambda+1)}{\alpha! \Gamma(\nb+\lambda+1)}}z^{\alpha}.
\end{equation*}
We will use also the weighted Bergman spaces $\mathcal{A}_{\lambda}^2(\mathbb{B}^{\nb'})$ and $\mathcal{A}_{\lambda+p+\nb'}^2(\mathbb{B}^{\nb''})$, with $p \in \mathbb{Z}_+$, whose basis elements are denoted by $e^{\lambda}_{\alpha'}(z')$ and $e^{\lambda+p+\nb'}_{\alpha''}(z'')$, respectively.
For each multi-index $\b=(\b_1, \cdots, \b_{m'})\in \mathbb{Z}_+^{m'}$ we introduce the Hilbert space
\begin{equation*}\label{GL_H_rho}
H_{\b}:= \overline{\textup{span}}\left\{ e_{\alpha}^{\lambda}(z) \, : \, \alpha=(\alpha_{(1)}, \cdots, \alpha_{(m')}, \alpha^{\prime \prime})
\in \mathbb{Z}_+^n \ \ \text{and} \ \ |\alpha_{(j)}|=\b_j, \ \ j =1,...,m'  \right\}.
\end{equation*}
Then we have the following orthogonal decomposition
\begin{equation*}\label{Orthogonal_decomposition_of_the_Bergman_space}
\mathcal{A}_{\lambda}^2(\mathbb{B}^\nb)=\bigoplus_{\b \in \mathbb{Z}_+^{m'}} H_{\b}.
\end{equation*}
A similar orthogonal decomposition can be performed for the Bergman space $\mathcal{A}^2_{\lambda}(\mathbb{B}^{\nb'})$:
\begin{equation*}
 \mathcal{A}^2_{\lambda}(\mathbb{B}^{\nb'}) = \bigoplus_{\b \in \mathbb{Z}^{m'}_+} \mathscr{H}_{\b},
\end{equation*}
where, for each $\b = (\b_1,...,\b_{m'}) \in \mathbb{Z}^{m'}_+$, the finite dimensional space $\mathscr{H}_{\b}$ is defined as
\begin{equation*}\label{defn_mathcal_H}
 \mathscr{H}_{\b} = \textup{span}\left\{ e_{\alpha^{\prime}}^{\lambda}(z') \, : \, \alpha^{\prime}=(\alpha_{(1)}, \cdots, \alpha_{(m')}) \in \mathbb{Z}_+^{\nb'} \ \ \text{and} \ \ |\alpha_{(j)}|=\b_j, \ \forall \: j =1,...,m'  \right\}.
\end{equation*}
Following \cite[Formula (2.6)]{BHaV}, for a multi-index $\b = (\b_1,...,\b_{m'}) \in \mathbb{Z}^{m'}_+$, we introduce the isometric isomorphism
\begin{equation*}
 \ub_{\b} \, : \ H_{\b} \ \ \longrightarrow \ \ \mathscr{H}_{\b} \otimes \mathcal{A}^2_{\lambda+|\b|+\ell}(\mathbb{B}^{\nb''}),
\end{equation*}
defined on the basis elements $e^{\lambda}_{\alpha}$, $\alpha = (\alpha^{\prime}, \alpha^{\prime\prime}) \in \mathbb{Z}^{\nb'}_+ \times \mathbb{Z}^{\nb''}_+$, of $H_{\b}$ by
\begin{equation*} \label{eq:u_rho}
 \ub_{\b} \, : \ e^{\lambda}_{\alpha}(z) \ \ \longmapsto \ \ e^{\lambda}_{\alpha^{\prime}}(z') \otimes e^{\lambda+|\b|+\nb'}_{\alpha^{\prime\prime}}(z'').
\end{equation*}
Next we introduce the Hilbert space
\begin{equation*} \label{eq:caligraphic_H}
 \mathcal{H} = \bigoplus_{\b \in \mathbb{Z}^{m'}_+} \mathcal{H}_{\b},\quad \text{with} \quad \mathcal{H}_{\b} = \mathscr{H}_{\b} \otimes \mathcal{A}^2_{\lambda+|\b|+\nb'}(\mathbb{B}^{\nb''})
\end{equation*}
and the unitary operator
\begin{equation*} \label{eq:U}
 U = \bigoplus_{\b \in \mathbb{Z}^{m'}_+} u_{\b} \, : \ \mathcal{A}_{\lambda}^2(\mathbb{B}^n)=\bigoplus_{\b \in \mathbb{Z}_+^{m'}} H_{\b} \ \longrightarrow \ \mathcal{H} = \bigoplus_{\b \in \mathbb{Z}^{m'}_+} \mathscr{H}_{\b} \otimes \mathcal{A}^2_{\lambda+|\b|+\nb'}(\mathbb{B}^{\nb''}),
\end{equation*}
acting component-wise with respect  to the direct sum decomposition.

Then Proposition 2.1 in \cite{BHaV} reads in our notations as follows.

\begin{proposition}[{\cite[Proposition 2.1]{BHaV}}]\label{Proposition_decomposition_Bergman_space}
 The unitary operator $U$ generates an isometric isomorphism between the spaces
\begin{equation*}
 \mathcal{A}_{\lambda}^2(\mathbb{B}^{\nb})=\bigoplus_{\b \in \mathbb{Z}_+^{m'}} H_{\b}
\quad {\rm and} \quad \mathcal{H} = \bigoplus_{\b \in \mathbb{Z}^{m'}_+} \mathscr{H}_{\b} \otimes \mathcal{A}^2_{\lambda+|\b|+\nb'}(\mathbb{B}^{\nb''}).
\end{equation*}
\end{proposition}

We return now to the symbols (in notations of \cite[Section 3]{BHaV}) of Toeplitz operators to be considered. Given a function $a \in L^{\infty}(\mathbb{B}^{\nb'})$ we denote by $f_a$ the function
\begin{equation*} \label{eq:f_a}
 a_w(z)\equiv f_a(z) \equiv f_a(z^{\prime},z^{\prime\prime}) = a\left(\frac{z^{\prime}}{\sqrt{1-|z^{\prime\prime}|^2}}\right) \in L^{\infty}(\mathbb{B}^\nb).
\end{equation*}
Similarly, for a function $b \in L^{\infty}(\mathbb{B}^{\nb''})$, we define the function $f_b(z) = f_b(z^{\prime},z^{\prime\prime}) = b(z^{\prime\prime}) \in L^{\infty}(\mathbb{B}^\nb)$.
In this notation, for the Toeplitz operators $\Tb^{\lambda}_{f_a}$ and $\Tb^{\lambda}_{f_c}$ we have, respectively,
\begin{eqnarray*}
 && \gamma_{f_a,\kn',\lambda}(\alpha) = \gamma_{a_w,\kn',\lambda}(\alpha') = \frac{\Gamma(\nb'+|\alpha'| + \lambda + 1)}{\Gamma(\lambda+1)  \prod_{j=1}^{m'} (k'_j-1+|\alpha'_{(j)}|)!} \\
 && \times \ \int_{\pmb{\Delta}_{m'}} a(\sqrt{\u_1},...,\sqrt{\u_{m'}}) (1-\u_1-...-\u_{m'})^\lambda
 \prod^{m'}_{j=1} {\u_j}^{|\alpha'_{(j)}|+k'_j-1} d\u_j,\\
 && \gamma_{f_b,\kn'',\lambda}(\alpha) = \gamma_{b,\kn'',\lambda}(\alpha',\alpha'') = \frac{\Gamma(\nb+|\alpha| + \lambda + 1)}{\Gamma(\l+|\a'|+\nb'+1)  \prod_{j=1}^{m''} (k''_j-1+|\alpha''_{(j)}|)!} \\
 && \times \ \int_{\tau(\Delta_{m''})}a(\sqrt{\r_1},...,\sqrt{\r_{m''}}) (1-\r_1...-\r_{m''})^{\l+ |\a'|+ \nb'} \prod^{m''}_{j=1} {\r_j}^{|\alpha''_{(j)}|+k''_j-1} d\r_j.
\end{eqnarray*}
Corollary 3.5 in \cite{BHaV} implies then  the following assertion.
\begin{proposition}[{\cite[Corollary 3.5]{BHaV}}]
The unitary operator $U$ realizes the unitary equivalence of Toeplitz operators $\Tb^{\lambda}_{f_{a}}$ and $\Tb^{\lambda}_{f_{b}}$ acting on $\mathcal{A}_{\lambda}^2(\mathbb{B}^\nb)$ with the following operators acting on $\mathcal{H} = \bigoplus_{\b \in \mathbb{Z}^{m'}_+} \mathscr{H}_{\b} \otimes \mathcal{A}^2_{\lambda+|\b|+\nb'}(\mathbb{B}^{\nb''})$:
\begin{eqnarray*}
 U \Tb^{\lambda}_{f_{a}} U^{-1} &=& \bigoplus_{\b \in \mathbb{Z}^{m'}_+}\,  \Tb_a^{\lambda}|_{\mathscr{H}_{\b}} \, \otimes \, I, \\
 U \Tb^{\lambda}_{f_{b}} U^{-1} &=& \bigoplus_{\b \in \mathbb{Z}^{m'}_+}\,  I \, \otimes \, \Tb_b^{\lambda+|\b|+\nb'}.
\end{eqnarray*}
\end{proposition}
Observe now that, for each  fixed $\a'=(\alpha'_{(1)},...,\alpha'_{(m')})$,
\begin{equation*}
 \Tb_a^{\lambda}|_{\mathscr{H}_{\a'}} = \gamma_{a_w,\kn',\lambda}(|\alpha'_{(1)}|,...,|\alpha'_{(m')}|)I,
\end{equation*}
so that $\Tb_a^{\lambda}|_{\mathscr{H}_{\a'}} = \gamma_{a_w,\kn',\lambda}(\alpha')I=\gamma_{a_w,\kn',\lambda}(|\alpha'_{(1)}|,...,|\alpha'_{(m')}|)$ is a scalar operator. On the other hand, the Toeplitz operator $ \Tb_b^{\lambda+|\b|+\nb'}$, acting on $\mathcal{A}^2_{\lambda+|\b|+\nb'}(\mathbb{B}^{\nb''})$, is unitary equivalent
with the multiplication operator $\gamma_{b,\kn'',\lambda+|\b|+\nb'}I$, acting on $\ell_2(\mathbb{Z}^{m''+1}_+)$. The entries in the sequence $\gamma_{b,\kn'',\lambda+|\b|+\nb'}$ (as it follows from \eqref{eq:gamma_b}) are given by
\begin{equation*}
 \gamma_{b,\kn'',\lambda+|\b|+\nb'}(\alpha'') = \gamma_{f_b,\kn'',\lambda}(\alpha) = \gamma_{b,\kn'',\lambda}(|\alpha'|,|\alpha''_{(1)}|,...,|\alpha''_{(m'')}|).
\end{equation*}
This explains the difference in the structure of the operators $\Tb^{\lambda}_{f_{a}}$ and $\Tb^{\lambda}_{f_{b}}$, which correspond to the above two cases, the special properties of the corresponding coefficients $\gamma$'s, and, as a consequence, the differences in the spectral representations of the Toeplitz operators in question.

\end{document}